\documentclass[12pt]{amsart}
\usepackage{amssymb,amsmath}

\numberwithin{equation}{section}

\newtheorem{Thm}{Theorem}[section]

\newtheorem{Lem}{Lemma}[section]

\def \eps{{\epsilon}}
\def \vareps{{\varepsilon}}

\def \sig{{\sigma}}

\def \xiHat{{\widehat{\xi}}}

\def \etaTld{{\tilde \eta}}
\def \thetaTld{{\tilde \theta}}
\def \rhoTld{{\tilde \rho}}
\def \chiTld{{\tilde \chi}}

\def \nuBar{{\overline{\nu}}}

\def \AA{{\mathcal A}}

\def \DD{{\mathcal D}}

\def \PP{{\mathcal P}}
\def \QQ{{\mathcal Q}}

\def \UU{{\mathcal U}}
\def \VV{{\mathcal V}}
\def \WW{{\mathcal W}}

\def \HBbb{{\mathbb H}}

\def \VBbb{{\mathbb V}}

\def \VHat{{\widehat V}}
\def \WHat{{\widehat W}}

\def \UTld{{\widetilde U}}
\def \VTld{{\widetilde V}}
\def \WTld{{\widetilde W}}

\def \YTld{{\widetilde Y}}

\def \pTld{{\tilde p}}

\def \uTld{{\tilde u}}
\def \vTld{{\tilde v}}

\def \BBar{{\overline{B}}}
\def \DBar{{\overline{D}}}
\def \QBar{{\overline{Q}}}

\def \DDBar{{\overline{\DD}}}
\def \QQBar{{\overline{\QQ}}}

\def \DDBarHat{{\widehat{\overline{\DD}}}}

\def \Dom{{\mathrm{Dom}}}
\def \Herm{{\mathrm{Herm}}}
\def \Null{{\mathrm{Null}}}

\def \Sp{{\mathrm{Sp}}}

\def \dee{{\mathrm d}}

\def \dt{{\dee t}}
\def \dx{{\dee x}}

\def \del{{\partial}}
\def \GRAD{\nabla_{\!\!x}}
\def \ROT{{\nabla_{\!\!x} \wedge}}

\def \DIV{\nabla_{\!\!x} \! \cdot }
\def \LAP{\Delta_x}

\def \DOT{{\,\cdot\,}}
\def \DDOT{{\,:\,}}

\def \To{{\, \longrightarrow \,}}

\def \<{\langle}
\def \>{\rangle}

\def\C{{\mathbb C}}      % complex numbers
\def\H{{\mathbb H}}      % Hilbert space
\def\N{{\mathbb N}}      % natural numbers
\def\R{{\mathbb R}}      % real numbers
\def\T{{\mathbb T}}      % torus
\def\Z{{\mathbb Z}}      % integers

\def\D{{\mathrm{D}}}

\def\RD{{{\mathbb R}^{\D}}}

\def\Ld{{{\mathbb L}^d}}

\def\Rd{{{\mathbb R}^d}}
\def\Sd{{{\mathbb S}^{d-1}}}
\def\Td{{{\mathbb T}^d}}
\def\Zd{{{\mathbb Z}^d}}

\def\NDim{{\mathrm{N}}}
\def\CN{{\C^{\NDim}}}
\def\CNN{{\C^{\NDim{\times}\NDim}}}
\def\RN{{\R^{\NDim}}}
\def\RNN{{\R^{\NDim{\times}\NDim}}}

\def\RDN{{\R^{d{\times}\NDim}}}
\def\RDNN{{\R^{d{\times}\NDim{\times}\NDim}}}
\def\RDDNN{{\R^{d{\times}d{\times}\NDim{\times}\NDim}}}

\def \H{{\mathbb H}}
\def \wH{\hbox{w-}\H}
\def \V{{\mathbb V}}
\def \wV{\hbox{w-}\V}

\def \kB{{\mathrm{k_B}}}  % Boltzmann Constant  --- R/N_A = 1.38065 10^{-23} J/K
\def \init{{\mathrm{in}}}
\def \loc{{\mathrm{loc}}}

\begin{document}

\title[Weakly Nonlinear-Dissipative Approximations]
      {Weakly Nonlinear-Dissipative Approximations
       of Hyperbolic-Parabolic Systems with Entropy}

\author[N. Jiang]%
       {Ning Jiang}
\address{Department of Mathematics \&
         Center for Scientific Computation and Mathematical Modeling
         (CSCAMM), University of Maryland, College Park, MD 20742}
\curraddr{Courant Institute of Mathematical Sciences,
          251 Mercer Street, New York, NY 10012, USA}
\email{njiang@cims.nyu.edu}

\author[C. D. Levermore]%
       {C. David Levermore}
\address{Department of Mathematics \&
         Institute for Physical Science and Technology (IPST),
         University of Maryland, College Park, MD 20742-4015, USA}
\email{lvrmr@math.umd.edu}

%\thanks{to be submitted to
%        {\em Archive of Rational Mechanics and Analysis}}

\date{\today}

\begin{abstract}

   Hyperbolic-parabolic systems have spatially homogenous stationary
states.  When the dissipation is weak, one can derive weakly
nonlinear-dissipative approximations that govern perturbations of
these constant states.  These approximations are quadratically
nonlinear.  When the original system has an entropy, the
approximation is formally dissipative in a natural Hilbert space.
We show that when the approximation is strictly dissipative it has
global weak solutions for all initial data in that Hilbert space.
We also prove a weak-strong uniqueness theorem for it.  In addition,
we give a Kawashima type criterion for this approximation to be
strictly dissipative.  We apply the theory to the compressible
Navier-Stokes system.

\end{abstract}

\maketitle

%%%%%%%%%%%%%%%%%%%%%%%%%%%%%%%%%%%%%%%%%%%%%%%%%%%%%%%%%%%%%%%%%%%%%%
%%  Section 1: Introduction                                         %%
%%%%%%%%%%%%%%%%%%%%%%%%%%%%%%%%%%%%%%%%%%%%%%%%%%%%%%%%%%%%%%%%%%%%%%

\section{Introduction}
\label{Introduction}

   We consider hyperbolic-parabolic systems over the $2\pi$-periodic
domain $\Td$ that have the form
\begin{equation}
  \label{Int-1}
  \del_t U + \DIV F(U) = \DIV \big[ D(U) \DOT \GRAD U \big] \,.
\end{equation}
These have spatially homogenous stationary solutions.  When the
dissipation is weak, one can derive a weakly nonlinear-dissipative
approximation that governs perturbations $\UTld$ about any constant
solution $U_o$.  These approximations have the form
\begin{equation}
  \label{Int-2}
  \del_t \UTld + \AA \UTld + \QQBar(\UTld,\UTld) = \DDBar \UTld \,,
\end{equation}
where $\AA=F_U(U_o)\DOT\GRAD$ is the linearization of the convection
operator about $U_o$, while quadratic operator $\QQBar$ and the linear
operator $\DDBar$ are formally given by
\begin{equation}
  \label{Int-3}
\begin{aligned}
  \QQBar(Y,Y)
  & = \lim_{T\to\infty}
      \frac{1}{2T}
      \int_{-T}^T e^{t \AA}
          \DIV \big[ \tfrac12 F_{UU}(U_o)(e^{-t \AA} Y, e^{-t \AA} Y)
                     \big] \, \dt \,,
\\
  \DDBar Y
  & = \lim_{T\to\infty}
      \frac{1}{2T}
      \int_{-T}^T e^{t \AA}
           D(U_o) \DDOT \GRAD^2 (e^{-t \AA} Y) \, \dt \,.
\end{aligned}
\end{equation}
Such weakly nonlinear-dissipative approximations arise when studying
incompressible limits of the compressible Navier-Stokes system
\cite{Danchin, Masmoudi}, global regularity of fast rotating
Navier-Stokes and Euler equations \cite{BMN97, BMN99, BMN01}, asymptotic
limits in equations of geophysical fluid dynamics \cite{EM1, EM2}, and
fast singular limits of hyperbolic and parabolic PDE's
\cite{Sh94, Gallagher1, Gallagher2}.

   We show that if the original system \eqref{Int-1} has a thrice
differentiable convex entropy structure then the approximating
system \eqref{Int-2} is formally dissipative in the Hilbert space
$\H$ whose inner product is given by
\begin{equation}
  \nonumber
  \big( \UTld \,|\, \VTld \big)_\H
  = \frac{1}{(2 \pi)^d}
    \int_\Td H_{UU}(U_o)\big( \UTld(x), \VTld(x) \big) \, \dx \,,
\end{equation}
where $H_{UU}(U_o)$ is the Hessian of the strictly convex entropy
density $H(U)$ at $U_o$.  This dissipation property follows because
the entropy structure implies that $\AA$ is skew-adjoint in $\H$,
that $\DDBar$ is nonpositive definite in $\H$, and that $\QQBar$
formally satisfies the cyclic identity
\begin{equation}
  \label{Int-4}
  0 = \big( \UTld \,|\, \QQBar\big( \VTld, \WTld \big) \big)_\H
      + \big( \VTld \,|\, \QQBar\big( \WTld, \UTld \big) \big)_\H
      + \big( \WTld \,|\, \QQBar\big( \UTld, \VTld \big) \big)_\H \,.
\end{equation}
We show that if $\DDBar$ is also strictly dissipative then the
approximating system \eqref{Int-2} has a Leray-type global weak
solution for all initial data in $\H$.  We cannot establish the
uniqueness of these solutions.  Indeed, when \eqref{Int-1} is the
Navier-Stokes system of gas dynamics then \eqref{Int-2} includes
the incompressible Navier-Stokes system as a subsystem.   The
uniqueness question therefore includes the uniqueness question
for the incompressible Navier-Stokes system.  We can however use
the cyclic identity \eqref{Int-4} to prove a so-called weak-strong
uniqueness theorem for the approximating system \eqref{Int-2}.

   It will be easily seen that $\DDBar$ will be nonnegative
definite if and only if the linear operators $\AA=F_U(U_o)\DOT\GRAD$
and $\DD=D(U_o)\DDOT\GRAD^2$ satisfy the Kawashima condition:
\begin{equation}
  \label{Int-5}
  \hbox{no nonconstant eigenfunction of $\AA$ }
  \hbox{is in the null space of $\DD$}.
\end{equation}
In one spatial dimension it is known that the Kawashima condition
implies that $\DDBar$ is strictly dissipative \cite{Kawashima}.
We give a stronger Kawashima-type criterion for $\DDBar$ to be
strictly dissipative in higher dimensions.

   Our paper is laid out as follows.  Section 2 presents the
entropy structure we will impose on system \eqref{Int-1}.  Section 3
presents the weakly nonlinear-dissipative approximation \eqref{Int-2}
and its relation to the Kawashima condition \eqref{Int-5}.  Section 4
develops the properties of the averaged operators \eqref{Int-3} that
we will need later.  This includes a proof of the cyclic identity
\eqref{Int-4}.  Section 5 presents our Kawashima-type criterion for
$\DDBar$ to be strictly dissipative.  Section 6 contains our existence
and weak-strong uniqueness theorems.  Finally, section 7 applies the
theory to the compressible Navier-Stokes system of gas dynamics.

%%%%%%%%%%%%%%%%%%%%%%%%%%%%%%%%%%%%%%%%%%%%%%%%%%%%%%%%%%%%%%%%%%%%%%
%%  Section 2: Hyperbolic-Parabolic Systems with Entropy            %%
%%%%%%%%%%%%%%%%%%%%%%%%%%%%%%%%%%%%%%%%%%%%%%%%%%%%%%%%%%%%%%%%%%%%%%

\section{Hyperbolic-Parabolic Systems with Entropy}
\label{HP}

   We consider hyperbolic-parabolic systems over $\Td$ in the
divergence form
\begin{equation}
  \label{HP-1}
  \del_t U + \DIV F(U) = \DIV \big[ D(U) \DOT \GRAD U \big] \,,
\end{equation}
where $U(x,t)$ is a density vector over $(x,t)\in\Td\times\R_+$ that
takes values in $\UU^c\subset\RN$.  Here $\Td=\Rd/(2\pi\Z)^d$ is the
$2\pi$-periodic torus and $\UU^c$ is the closure of a convex domain
$\UU\subset\RN$.  We assume that the flux tensor $F:\UU\to\RDN$ is
twice continuously differentiable such that
\begin{equation}
  \label{HP-2}
  \del_t U + \DIV F(U) = 0 \quad
  \hbox{is hyperbolic} \,,
\end{equation}
while the diffusion tensor $D:\UU\to\RDDNN$ is continuously
differentiable such that
\begin{equation}
  \label{HP-3}
  \del_t U = \DIV [D(U) \DOT \GRAD U] \quad
  \hbox{is parabolic} \,.
\end{equation}
Recall that system (\ref{HP-2}) is said to be hyperbolic if for
every $U\in\UU$ and every $\xi\in\Rd$ the $N{\times}N$ matrix
$F_U(U)\DOT\xi$ is diagonalizable within the reals --- i.e. it has a
complete set of real eigenvectors.  System (\ref{HP-2}) is said to
be strictly hyperbolic if moreover the eigenvalues of
$F_U(U)\DOT\xi$ are distinct.   Recall that system (\ref{HP-3}) is
said to be parabolic if for every $U\in\UU$ and every $\xi\in\Rd$
the $N{\times}N$ matrix $D(U)\DDOT\xi^{\otimes2}$ is diagonalizable
within the reals and has nonnegative eigenvalues.  System
(\ref{HP-3}) is said to be strictly parabolic if moreover the
eigenvalues of $D(U)\DDOT\xi^{\otimes2}$ are positive.  Many studies
of hyperbolic-parabolic systems assume that system (\ref{HP-2}) is
strictly hyperbolic while system (\ref{HP-3}) is strictly parabolic.
We will not do that here.  Rather, we will assume that system
(\ref{HP-1}) has a strictly convex entropy and satisfies certain
nonsingularity conditions.

\subsection{Entropy Structure}
\label{HP-Entropy}
We say that $H:\UU\to\R$ is a strictly convex entropy for the system
(\ref{HP-1}) when $H$ is twice continuously differentiable over $\UU$
and for every $U\in\UU$
\begin{equation}
  \label{HP-4}
\begin{aligned}
  \hbox{(i)} \,
  & H_{UU}(U) &
  & \hbox{is positive definite} \,,
\\
  \hbox{(ii)} \,
  & H_{UU}(U) \, F_U(U) \DOT \xi &
  & \hbox{is symmetric for every $\xi\in\Rd$} \,,
\\
  \hbox{(iii)} \,
  & H_{UU}(U) \, D(U) &
  & \hbox{is symmetric and nonnegative definite} \,.
\end{aligned}
\end{equation}
The existence of such a strictly convex entropy implies that system
(\ref{HP-1}) is hyperbolic-parabolic.  The compressible
Navier-Stokes system is a hyperbolic-parabolic system that is
neither strictly hyperbolic nor strictly parabolic, yet has a
strictly convex entropy. We will study this example in Section
\ref{CNS}.   There are many other systems from physics that fit into
this framework \cite{CLL}.

   Condition (ii) in (\ref{HP-4}) implies there exists $J:\UU\to\Rd$
that is twice continuously differentiable such that
\begin{equation}
  \label{HP-5}
  H_U(U) F_U(U) = J_U(U) \,.
\end{equation}
It follows that if $U(x,t)$ is a classical solution of (\ref{HP-1})
that takes its values in $\UU$ then it satisfies
\begin{equation}
  \label{HP-6}
  \del_t H(U) + \DIV J(U)
  = \DIV \big[ H_U(U) \, D(U) \DOT \GRAD U \big]
    - \GRAD H_U(U) \DOT D(U) \DOT \GRAD U \,.
\end{equation}

   Condition (iii) in (\ref{HP-4}) implies that for every
differentiable $U:\Td\to\RN$ with $U(x)\in\UU$
\begin{equation}
  \label{HP-7}
  \GRAD H_U(U) \DOT D(U) \DOT \GRAD U
  = \GRAD U^T \DOT  H_{UU}(U) D(U) \DOT \GRAD U \geq 0 \,.
\end{equation}
We thereby see that (\ref{HP-6}) is a local dissipation law for
$H(U)$.  When (\ref{HP-6}) is integrated over $\Td$ one obtains
the global dissipation law
\begin{equation}
  \label{HP-8}
\begin{aligned}
  \frac{\dee}{\dt} \int H(U) \, \dx
  & = - \int \GRAD H_U(U) \DOT D(U) \DOT \GRAD U \, \dx
\\
  & = - \int \GRAD U^T \DOT  H_{UU}(U) D(U) \DOT \GRAD U \, \dx
    \leq 0 \,.
\end{aligned}
\end{equation}

\subsection{Nonsingularity Condition}
\label{HP-Nonsingularity}

   The class of stationary classical solutions of (\ref{HP-1})
that take values in $\UU$ is constrained by the entropy structure.
It follows from (\ref{HP-8}) that every such solution satisfies
\begin{equation}
  \nonumber
  \int \GRAD U^T \DOT  H_{UU}(U) D(U) \DOT \GRAD U \, \dx = 0 \,,
\end{equation}
which by (\ref{HP-7}) implies that
\begin{equation}
  \nonumber
  \GRAD U^T \DOT  H_{UU}(U) D(U) \DOT \GRAD U = 0 \,.
\end{equation}
Conditions (i) and (iii) then imply that $D(U)\DOT\GRAD U=0$,
which when plugged into (\ref{HP-1}) with $\del_tU=0$ yields
$\DIV F(U)=0$.  We thereby see that every classical solution of
(\ref{HP-1}) that takes its values in $\UU$ satisfies
\begin{equation}
  \label{HP-9}
  \DIV F(U) = 0 \,, \qquad
  D(U) \DOT \GRAD U = 0 \,.
\end{equation}
In general this is not enough information to conclude that $U$
is a constant.

   We consider hyperbolic-parabolic systems over $\Td$ that are
nonsingular in the sense that for every $U\in\UU$ and every
continuously differentiable $\UTld:\Td\to\RN$ one has that
\begin{equation}
  \label{HP-10}
\left.
\begin{aligned}
  F_U(U) \DOT \GRAD \UTld & = 0
\\
  D(U) \DOT \GRAD \UTld & = 0
\end{aligned} \,
\right\} \quad \implies \quad
  \GRAD \UTld = 0 \,.
\end{equation}
With this additional assumption we can prove the following.

\begin{Lem}
\label{HP-Lem-1}
If the hyperbolic-parabolic system (\ref{HP-1}) has a strictly convex
entropy (\ref{HP-4}) and satisfies the nonsingularity condition
(\ref{HP-10}) then its only stationary, spatially periodic classical
solutions that take their values in $\UU$ are constant solutions.
\end{Lem}

\noindent In the next section we will strengthen the nonsingularity
condition (\ref{HP-10}).

%%%%%%%%%%%%%%%%%%%%%%%%%%%%%%%%%%%%%%%%%%%%%%%%%%%%%%%%%%%%%%%%%%%%%%
%%  Section 2: Weakly Nonlinear-Dissipative Approximations          %%
%%%%%%%%%%%%%%%%%%%%%%%%%%%%%%%%%%%%%%%%%%%%%%%%%%%%%%%%%%%%%%%%%%%%%%

\section{Weakly Nonlinear-Dissipative Approximations}
\label{WND}

   We now consider
{\em weakly nonlinear-dissipative (WND) approximations} to
hyperbolic-parabolic systems of the form (\ref{HP-1}).  In order to
see how these approximations depend on the choice of dependent
variables, we express (\ref{HP-1}) in the form
\begin{equation}
  \label{WND-1}
  \del_t U(W) + \DIV F(U(W))
  = \DIV \big[ D(U(W)) \DOT \GRAD U(W) \big] \,.
\end{equation}
Here the mapping $U:\WW\to\UU$ is assumed to be a twice continuously
differentiable bijection and have a nonsingular Jacobian.  This means
that its inverse mapping $U^{-1}:\UU\to\WW$ will also be a twice
continuously differentiable bijection and have a nonsingular Jacobian.

\subsection{Linearization and the Kawashima Condition}

   In order to motivate our approximation, let us first consider the
linearization of the hyperbolic-parabolic system (\ref{WND-1}) about
some constant state $W_o\in\WW$.  This is
\begin{equation}
  \label{WND-2}
  \del_t \WTld + A_o \DOT \GRAD \WTld = B_o \DDOT \GRAD^2 \WTld \,,
\end{equation}
where $A_o\in\RDNN$, and $B_o\in\RDDNN$ are defined by
\begin{equation}
  \label{WND-3}
  A_o = R_o^{-1} F_U(U_o) R_o \,, \qquad
  B_o = R_o^{-1} D(U_o) R_o \,,
\end{equation}
with $U_o\in\UU$ and $R_o\in\RNN$ given by
\begin{equation}
  \label{WND-4}
  U_o = U(W_o) \,, \qquad
  R_o = \del_V U(W_o) \,.
\end{equation}

   If $H:\UU\to\R$ is a strictly convex entropy for system
(\ref{HP-1}) then system (\ref{WND-2}) is symmetrized by the
positive definite matrix
\begin{equation}
  \label{WND-5}
  G_o = R_o^T H_{UU}(U_o) R_o \,.
\end{equation}
Specifically, one sees from (\ref{HP-4}) and (\ref{WND-3}) that
\begin{equation*}
\begin{aligned}
  G_o A_o \DOT \xi & = R_o^T H_{UU}(U_o) F_U(U_o) \DOT \xi R_o & \quad
  & \hbox{is symmetric for every $\xi\in\Rd$} \,,
\\
   G_o B_o & = R_o^T H_{UU}(U_o) D(U_o) R_o & \quad
  & \hbox{is symmetric and nonnegative definite} \,.
\end{aligned}
\end{equation*}
The solutions of (\ref{WND-2}) thereby satisfy the local dissipation
law
\begin{equation*}
  \del_t \big( \tfrac12 \WTld^T G_o \WTld \big)
  + \DIV \big( \tfrac12 \WTld^T G_o A_o \WTld \big)
  = \DIV \big( \WTld^T G_o B_o \DOT \GRAD \WTld \big)
    - \GRAD \WTld^T \DOT G_o B_o \DOT \GRAD \WTld \,.
\end{equation*}
When this equation is integrated over $\Td$ one obtains the global
dissipation law
\begin{equation*}
  \frac{\dee}{\dt} \int \tfrac12 \WTld^T G_o \WTld \, \dx
  = - \int \GRAD \WTld^T \DOT G_o B_o \DOT \GRAD \WTld \, \dx
  \leq 0 \,.
\end{equation*}
The initial-value problem for system (\ref{WND-2}) is therefore
naturally well-posed in the Hilbert space $\H=L^2(\dx;\RN)$
equipped with the inner product
\begin{equation}
  \nonumber
  \big( \WTld_1 \,|\, \WTld_2 \big)_\HBbb
  = \frac{1}{(2 \pi)^d} \int_\Td \WTld_1^T G_o \WTld_2 \, \dx \,,
\end{equation}
where $G_o$ is given by (\ref{WND-5}).

  By arguing as we did in the nonlinear settings, every stationary,
spatially periodic classical solution of system (\ref{WND-2}) must
satisfy
\begin{equation}
  \nonumber
  A_o \DOT \GRAD \WTld = 0 \,, \qquad
  B_o \DOT \GRAD \WTld = 0 \,.
\end{equation}
We can conclude that $\GRAD\WTld=0$ if system (\ref{HP-1}) satisfies
the nonsingularity condition (\ref{HP-10}).  In that case the only
stationary, spatially periodic classical solutions of system
(\ref{WND-2}) are the constant solutions.

   It is natural to ask if every solution of system (\ref{WND-2}) in
$\H$ will decay to a constant solution as $t\to\infty$.  Kawashima
\cite{Ka,KaSh} gave an elegant characterization of when this is the
case in terms of the skew-adjoint operator $\AA$ and the self-adjoint
operator $\DD$ that are formally given by
\begin{equation}
  \label{WND-6}
  \AA = A_o \DOT \GRAD \,, \qquad
  \DD = B_o \DDOT \GRAD^2 \,.
\end{equation}
Because $\AA=A_o\DOT\GRAD$ is skew-adjoint in $\H$, its spectrum,
$\Sp(\AA)$, is purely imaginary.  When the spatial domain is bounded
$\Sp(\AA)$ only contains eigenvalues.  Moreover, these eigenvalues
and their corresponding eigenfunctions are easily computed by
Fourier methods when the spatial domain is $\Td$.  Kawashima gave
the following characterization.

\begin{Thm}
\label{WND-Thm-1}
Every solution of system (\ref{WND-2}) in $\H$ will decay to a
constant solution as $t\to\infty$ if and only if
\begin{equation}
  \label{WND-7}
  \hbox{no nonconstant eigenfunction of $\AA$ }
  \hbox{is in the null space of $\DD$}.
\end{equation}
\end{Thm}

\smallskip
\noindent {\bf Remark.}  The Kawashima condition (\ref{WND-7}) is
clearly necessary for every solution of system (\ref{WND-2}) in $\H$
to decay to a constant solution as $t\to\infty$.  Indeed, if $V$ is
nonconstant eigenfunction of $\AA$ for the eigenvalue $i\omega$ such
that $V$ lies in the null space of $\DD$ then the real part of
$e^{-i\omega t}V$ is a real-valued solution of (\ref{WND-2}) that
does not decay to a constant solution as $t\to\infty$.

\smallskip
\noindent
{\bf Remark.}  The Kawashima condition (\ref{WND-7}) is stronger
than our nonsingular condition (\ref{HP-10}).  It has been used to
obtain similar results regarding the existence, regularity, and decay
as $t\to\infty$ of global solutions to nonlinear systems
\cite{BHN, HZ1, HZ2, Kawashima, Kawashima-BE, Ru-Serre, Yong-1}.
Villani has generalized it in his theory of hypercoercivity
\cite{Villani} where he discusses its relation to the older
H\"{o}rmander hypoellipticity condition \cite{Hormander}.

\subsection{Weak Dissipation, Two-Time Asymptotics}

   It is more interesting to consider regimes in which the
dissipation is weak.  Introduce the nondimensional (small) parameter
$\eps>0$ so that
\begin{equation}
  \label{WND-8}
  \del_t U(W_\eps) + \DIV F(U(W_\eps))
   = \eps \, \DIV \big[ D(U(W_\eps)) \DOT \GRAD U(W_\eps) \big] \,,
\end{equation}
Let $W_o\in\WW$ and set
\begin{equation*}
  W_\eps = W_o + \eps \, \WTld + \eps^2 \WTld^{(2)} + O(\eps^3) \,.
\end{equation*}
Then
\begin{equation*}
  U(W_\eps) = U_o + \eps \, R_o \WTld
              + \eps^2 R_o \big( \WTld^{(2)} + S_o(\WTld,\WTld) \big)
              + O(\eps^3) \,,
\end{equation*}
where $U_o$ and $R_o$ are given by (\ref{WND-4}) while $S_o$
is defined by
\begin{equation*}
  S_o(\WTld,\WTld)
  = \tfrac12 R_o^{-1} \del_{WW} U(W_o)(\WTld,\WTld) \,.
\end{equation*}

    To first order in $\eps$ we see that $\WTld$ satisfies
\begin{equation}
  \label{WND-9}
  \del_t \WTld + A_o \DOT \GRAD \WTld = 0 \,.
\end{equation}
The solution of this equation is given by
\begin{equation*}
  \WTld(t) = e^{-t\AA} \WTld^\init \,, \quad
  \hbox{where $\AA=A_o \DOT \GRAD$} \,.
\end{equation*}
Because $\AA$ is skew-adjoint on $\HBbb$, the solution operator
$e^{-t\AA}$ is strongly continuous, one parameter group of unitary
operators on $\HBbb$.  This approximation cannot be valid uniformly
in time because (1) the solutions of (\ref{WND-9}) do not decay as
$t\to\infty$ and (2) there are generally nonconstant stationary
solutions of (\ref{WND-9}) when $d\geq2$.

   In order to overcome these problems, one has to introduce a slow
time scale $\tau=\eps\,t$ into the asymptotics and consider
\begin{equation}
  \nonumber
  \eps \, \del_\tau U(W_\eps) + \del_t U(W_\eps) + \DIV F(U(W_\eps))
   = \eps \, \DIV \big[ D(U(W_\eps)) \DOT \GRAD U(W_\eps) \big] \,.
\end{equation}
To first order we see that $\WTld$ still satisfies (\ref{WND-8}).
Hence,
\begin{equation}
  \label{WND-10}
  \WTld(t,\tau) = e^{-t \AA} \YTld(\tau) \,,
\end{equation}
where the $\tau$ dependence of $\YTld$ has yet to be determined.

   To the second order in $\eps$ we see that
\begin{equation}
  \label{WND-11}
  \del_t \WTld^{(2)} + \AA \WTld^{(2)}
  = - \del_\tau \WTld - \DIV Q_o(\WTld,\WTld)
    + B_o \DDOT \GRAD^2 \WTld
    - (\del_t + \AA) S_o(\WTld,\WTld) \,,
\end{equation}
where
\begin{equation}
  \label{WND-12}
  Q_o(\WTld,\WTld)
  = \tfrac12 R_o^{-1} F_{UU}(U_o)(R_o \WTld, R_o \WTld) \,.
\end{equation}
The right-hand side of (\ref{WND-11}) is an almost periodic
function of $t$.  For $\WTld^{(2)}$ to be an almost periodic, we must
require
\begin{equation}
  \nonumber
  \lim_{T\to\infty} \frac1{2T} \int_{-T}^T
  e^{t\AA} \big[ \del_\tau \WTld + \DIV Q_o(\WTld,\WTld)
                 - B_o \DDOT \GRAD^2 \WTld \big] \, \dt = 0 \,.
\end{equation}
Hence, because $\WTld$ is given by (\ref{WND-10}), we see that
$\YTld(\tau)$ satisfies
\begin{equation}
  \label{WND-13}
  \del_\tau \YTld + \QQBar(\YTld,\YTld) = \DDBar \YTld \,,
\end{equation}
where the operators $\QQBar$ and $\DDBar$ are formally defined by
\begin{equation}
  \label{WND-14}
\begin{aligned}
  \QQBar(Y,Y)
  & = \lim_{T\to\infty}
      \frac{1}{2T} \int_{-T}^T e^{t \AA}
                       \QQ(e^{-t \AA} Y,e^{-t \AA} Y) \, \dt \,,
\\
  \DDBar Y
  & = \lim_{T\to\infty}
      \frac{1}{2T} \int_{-T}^T e^{t \AA}
                       \DD e^{-t \AA} Y \, \dt \,,
\end{aligned}
\end{equation}
with the operators $\QQ$ and $\DD$ being given by
\begin{equation*}
  \QQ(W,W) = \DIV Q_o(W,W) \,, \qquad
  \DD W = B_o \DDOT \GRAD^2 W \,.
\end{equation*}

   It is easily checked from formulas (\ref{WND-14}) that
\begin{equation}
  \label{WND-15}
  \QQBar(e^{-t\AA} W, e^{-t\AA} W)
  = e^{-t\AA} \QQBar(W,W) \,, \qquad
  \DDBar e^{-t\AA} = e^{-t\AA} \DDBar \,.
\end{equation}
It then follows from (\ref{WND-10}) and (\ref{WND-13}) that
$\WTld(t,\eps\,t)$ satisfies
\begin{equation*}
  \del_t \WTld + \AA \WTld + \eps \, \QQBar(\WTld,\WTld)
  = \eps \, \DDBar \WTld \,.
\end{equation*}
Setting $\eps=1$, we call this the weakly nonlinear-dissipative
approximation of the hyperbolic-parabolic system (\ref{WND-1}).

\subsection{Change of Dependent Variables}

   Suppose we had expressed the system (\ref{WND-6}) in terms of
different dependent variables $W'$:
\begin{equation}
  \label{WND-16}
  \del_t U'(W'_\eps) + \DIV F(U'(W'_\eps))
   = \eps \, \DIV \big[ D(U'(W'_\eps)) \DOT \GRAD U'(W'_\eps) \big] \,,
\end{equation}
where $U':\WW'\to\UU$ is a twice continuously differentiable bijection
with a nonsingular Jacobian.  Let $W'_o$ be the unique constant state
such that $U'(W'_o)=U_o=U(W_o)$.  If we approximate solutions of
(\ref{WND-16}) that are near $W'_o$ just as we approximated solutions
of (\ref{WND-8}) that are near $W_o$ then to leading order we obtain
\begin{equation*}
  \del_t \WTld' + A'_o \DOT \GRAD \WTld' = 0 \,,
\end{equation*}
where $A'_o = T_o^{-1} A_o T_o$ with
\begin{equation}
  \label{WND-17}
  T_o = \del_{W'} \big(U^{-1}\big(U'(W')\big)\big) \Big|_{W'=W'_o}
      = R_o^{-1} R'_o \,.
\end{equation}
In other words, the leading order approximation of the transformed
system (\ref{WND-16}) is the leading order approximation of the
original system transformed by the associated linear change of
variables $\WTld=T_o\WTld'$.

   Remarkably, the same transformation property holds for the
weakly nonlinear-dissipative approximation!  Specifically, if the
weakly nonlinear-dissipative approximation of the primed system
(\ref{WND-16}) is
\begin{equation*}
  \del_t \WTld' + \AA' \WTld' + \eps \, \QQBar'(\WTld',\WTld')
  = \eps \, \DDBar' \WTld \,,
\end{equation*}
then it is related to the weakly nonlinear-dissipative approximation
of the unprimed system by
\begin{equation}
  \label{WND-18}
  \AA' = T_o^{-1} \AA T_o \,, \qquad
  \QQBar'(\WTld',\WTld')
  = T_o^{-1} \QQBar(T_o \WTld', T_o \WTld') \,, \qquad
  \DDBar' = T_o^{-1} \DDBar T_o \,,
\end{equation}
where $T_o$ is again given by (\ref{WND-17}).  In other words, it is
simply the original nonlinear-dissipative approximation transformed by
the associated linear change of variables $\WTld=T_o\WTld'$.  This
fact allows us to derive the weakly nonlinear-dissipative approximation
of system (\ref{WND-8}) in any set of dependent variables we choose
becasue the result is unique up to the transformation (\ref{WND-18}).

\section{Averaged Operators}

   In this section we collect some properties of the averaged
operators $\QQBar$ and $\DDBar$.

\subsection{Spectral Formulas for the Averaged Operators}

   Because $\AA=A_o\DOT\GRAD$ is skew-adjoint in $\H$, its spectrum,
$\Sp(\AA)$, is purely imaginary.  Define
\begin{equation*}
  \sig(\AA)
  = \big\{ \omega \in \R \,:\, i \omega \in \Sp(\AA) \big\} \,.
\end{equation*}
For every $\omega\in\sig(\AA)$ let $\H_\omega$ denote the eigenspace
associated with the eigenvalue $i\omega$ of $\AA$ and let $E_\omega$
be the $\H$-orthogonal projection onto $\H_\omega$.  Then for every
$\omega\in\sig(\AA)$ one has
\begin{equation*}
\begin{gathered}
  (E_\omega W_1 \,|\, W_2)_\H = (W_1 \,|\, E_\omega W_2)_\H \quad
  \hbox{for every $W_1,W_2\in\H$} \,,
\\
  E_\omega^{\,2} = E_\omega \,, \qquad \qquad
  \H_{\omega} = E_\omega \H \,.
\end{gathered}
\end{equation*}
For every $\omega\in\sig(\AA)$ let $C_\omega^\infty$ denote the
smooth functions in $\H_\omega$ --- i.e. let
$C_\omega^\infty=\H_\omega\cap C^\infty(\Td)$.  We will assume that
each $C_\omega^\infty$ is dense in $\H_\omega$ and is contained in
the domains of $\DD$ and $\QQ$.  This assumption holds because we
are working over the periodic domain $\Td$.  In general settings it
would mean that every smooth eigenfunction of $\AA$ must also
satisfy any boundary conditions associated with $\DD$, which is
typically is not true.  Our periodic setting we have explicit
formulas for each $E_\omega$, however we will put off using them
as long as possible.

   The Spectral Decomposition Theorem implies that
\begin{equation*}
  e^{-t\AA} = \sum_{\omega\in\sig(\AA)}
                  e^{- i \omega t} E_\omega \quad
  \hbox{for every $t\in\R$} \,.
\end{equation*}
If $W\in\H$ has nonzero components only in a finite number of the
$\H_\omega$ then there are a finite number of nonzero terms in the
decomposition
\begin{equation*}
  e^{-t\AA} W = \sum_{\omega\in\sig(\AA)}
                    e^{- i \omega t} E_\omega W \,.
\end{equation*}
If we assume moreover that each of these components is smooth then
by using this decomposition in (\ref{WND-14}) we obtain the spectral
formulas
\begin{equation}
  \label{AOp-1}
\begin{aligned}
  \QQBar(W,W)
  & = \sum_{\omega_1,\,\omega_2\in\sig(\AA)}
          E_{\omega_1+\omega_2} \QQ(E_{\omega_1} W, E_{\omega_2} W) \,,
\\
  \DDBar W
  & = \sum_{\omega\in\sig(\AA)} E_\omega \DD E_\omega W \,,
\end{aligned}
\end{equation}
where we understand that $E_{\omega_1+\omega_2}=0$ when
$\omega_1+\omega_2\not\in\sig(\AA)$.

\subsection{Bounds on the Averaged Operators}

   The spectral formulas (\ref{AOp-1}) only apply to those $W\in\H$ that
have nonzero components in only a finite number of the $\H_\omega$,
each of which is smooth.  Denote this set by $\oplus_\omega
C_\omega^\infty$.  While this set is dense in $\H$, it must be
extended to larger classes of $W$.  This is done by continuity once
we obtain appropriate bounds on the forms associated with the
operators $\QQBar$ and $\DDBar$.

   We begin by recalling some bounds on the forms associated with the
operators $\AA$, $\QQ$, and $\DD$ in terms of norms that are
invariant under $e^{-t\AA}$.  Specifically, we employ the $\H^s$
spaces that are the completion of $\oplus_\omega C_\omega^\infty$ in
the norms defined for every $s\in\R$ by
\begin{equation}
  \nonumber
  \| W \|_{\H^s}
  = \left( \sum_{\xi\in\Zd} \big( 1 + |\xi|^2 \big)^s |\WHat(\xi)|_{G_o}
           \right)^\frac12 \,, \qquad
  \WHat(\xi) = \frac{1}{(2 \pi)^d}
               \int_\Td e^{-i\xi\DOT x} W(x) \, \dx \,.
\end{equation}
Here $\WHat(\xi)$ is the Fourier coefficient associated with the wave
vector $\xi\in\Zd$.

\begin{Lem}
\label{AOp-Lem-1}
There there exist positive constants $C_\AA$ and $C_\DD$ such that for
every $W_1,W_2\in C^\infty(\Td)$ and every $s'\in\R$ one has the bounds
\begin{equation}
  \label{AOp-2}
\begin{aligned}
  \big| \big( W_1 \,|\, \AA W_2 \big)_\H \big|
  & \leq C_\AA \| W_1 \|_{\H^{-s'}}
               \| \GRAD W_2 \|_{\H^{s'}} \,.
\\
  \big| \big( W_1 \,|\, \DD W_2 \big)_\H \big|
  & \leq C_\DD \| \GRAD W_1 \|_{\H^{-s'}}
               \| \GRAD W_2 \|_{\H^{s'}} \,.
\end{aligned}
\end{equation}
For every $s>d/2$ there exist a positive constant $C^s_\QQ$ such that
for every $W_1,W_2,W_3\in C^\infty(\Td)$ and every $s'\geq0$ one has
the bounds
\begin{equation}
  \label{AOp-3}
\begin{aligned}
  \big| \big( W_1 \,|\, \QQ(W_2,W_3) \big)_\H \big|
  & \leq C^s_\QQ \| \GRAD W_1 \|_{\H^s}
                 \| W_2 \|_\H \| W_3 \|_\H \,,
\\
  \big| \big( W_1 \,|\, \QQ(W_2,W_3) \big)_\H \big|
  & \leq 2^{s'} C^s_\QQ
         \| W_1 \|_{\H^{-s'}}
         \Big( \| \GRAD W_2 \|_{\H^{s'}} \| W_3 \|_{\H^s}
               + \| W_2 \|_{\H^s} \| \GRAD W_3 \|_{\H^{s'}} \Big) \,.
\end{aligned}
\end{equation}
\end{Lem}

\smallskip
\noindent
{\bf Proof.} We refer the reader to \cite{CF} for similar proofs. \qed

\smallskip

   Because the $\H^s$ norms are invariant under $e^{-t\AA}$ we can prove
the following.

\begin{Lem}
\label{AOp-Lem-2}
For every $s>d/2$, $s'\geq0$, and
$W_1,W_2,W_3\in\oplus_\omega C_\omega^\infty$ one has the bounds
\begin{equation}
   \label{AOp-4}
\begin{aligned}
  \big| \big( W_1 \,|\, \DDBar W_2 \big)_\H \big|
  & \leq C_\DD \| \GRAD W_1 \|_{\H^{-s'}}
               \| \GRAD W_2 \|_{\H^{s'}} \,,
\\
  \big| \big( W_1 \,|\, \QQBar(W_2,W_3) \big)_\H \big|
  & \leq C^s_\QQ \| \GRAD W_1 \|_{\H^s}
                 \| W_2 \|_\H \| W_3 \|_\H \,,
\\
  \big| \big( W_1 \,|\, \QQBar(W_2,W_3) \big)_\H \big|
  & \leq 2^{s'} C^s_\QQ
         \| W_1 \|_{\H^{-s'}}
         \Big( \| \GRAD W_2 \|_{\H^{s'}} \| W_3 \|_{\H^s}
               + \| W_2 \|_{\H^s} \| \GRAD W_3 \|_{\H^{s'}} \Big) \,,
\end{aligned}
\end{equation}
where $C_\DD$ and $C^s_\QQ$ are the constants appearing in
(\ref{AOp-2}) and (\ref{AOp-3}) of Lemma \ref{AOp-Lem-1}.
\end{Lem}

\smallskip
\noindent
{\bf Proof.} Because the $\H^s$ norms are invariant under $e^{-t\AA}$,
and because $e^{-t\AA}$ and $\GRAD$ commute, for every $s\geq0$ and
every $W\in C^\infty(\Td)$ we have
\begin{equation}
  \label{AOp-5}
  \big\| e^{-t \AA} W \big\|_{\H^s}
  = \big\| W \big\|_{\H^s} \,, \qquad
  \big\| \GRAD e^{-t \AA} W \big\|_{\H^s}
  = \big\| \GRAD W \big\|_{\H^s} \,.
\end{equation}

   We first prove the bound on $\DDBar$ in (\ref{AOp-4}).  From
the definition of $\DDBar$ given in (\ref{WND-14}) we see that
\begin{equation*}
\begin{aligned}
  \big( W_1 \,|\, \DDBar W_2 \big)_\H
  & = \lim_{T\to\infty}
      \frac{1}{2T} \int_{-T}^T
                       \big( W_1 \,|\,
                             e^{t \AA} \DD e^{-t \AA} W_2
                             \big)_\H \, \dt
\\
  & = \lim_{T\to\infty}
      \frac{1}{2T} \int_{-T}^T
                       \big( e^{-t \AA} W_1 \,|\,
                             \DD e^{-t \AA} W_2 \big)_\H \, \dt \,.
\end{aligned}
\end{equation*}
Hence, the bound on $\DD$ in (\ref{AOp-2}) and the invariances
(\ref{AOp-5}) imply
\begin{equation*}
\begin{aligned}
  \big| \big( W_1 \,|\, \DDBar W_2 \big)_\H \big|
  & \leq \lim_{T\to\infty}
         \frac{1}{2T} \int_{-T}^T
                          \big| \big( e^{-t \AA} W_1 \,|\,
                                      \DD e^{-t \AA} W_2
                                      \big)_\H \big| \, \dt
\\
  & \leq \lim_{T\to\infty}
         \frac{C_\DD}{2T}
         \int_{-T}^T
             \big\| \GRAD e^{-t \AA} W_1 \big\|_\H
             \big\| \GRAD e^{-t \AA} W_2 \big\|_\H \, \dt
\\
  & = \lim_{T\to\infty}
      \frac{C_\DD}{2T}
      \int_{-T}^T \| \GRAD W_1 \|_\H  \| \GRAD W_2 \|_\H \, \dt
    = C_\DD \| \GRAD W_1 \|_\H \| \GRAD W_2 \|_\H \,.
\end{aligned}
\end{equation*}
This proves the bound on $\DDBar$ in (\ref{AOp-4}).

   We now prove the bounds on $\QQBar$ in (\ref{AOp-4}).  From
the definition of $\QQBar$ given in (\ref{WND-14}) we see that
\begin{equation*}
\begin{aligned}
  \big( W_1 \,|\, \QQBar(W_2,W_3) \big)_\H
  & = \lim_{T\to\infty}
      \frac{1}{2T} \int_{-T}^T \big( W_1 \,|\,
                                     e^{t \AA} \QQ(e^{-t \AA} W_2,
                                                   e^{-t \AA} W_3)
                                     \big)_\H \, \dt
\\
  & = \lim_{T\to\infty}
      \frac{1}{2T} \int_{-T}^T \big( e^{-t \AA} W_1 \,|\,
                                     \QQ(e^{-t \AA} W_2,
                                         e^{-t \AA} W_3)
                                     \big)_\H \, \dt \,.
\end{aligned}
\end{equation*}
Hence, the first bound on $\QQ$ in (\ref{AOp-3}) and the invariances
(\ref{AOp-5}) imply
\begin{equation*}
\begin{aligned}
  \big| \big( W_1 \,|\, \QQBar(W_2,W_3) \big)_\H \big|
  & \leq \lim_{T\to\infty}
         \frac{1}{2T} \int_{-T}^T
                          \big| \big( e^{-t \AA} W_1 \,|\,
                                      \QQ(e^{-t \AA} W_2,
                                          e^{-t \AA} W_3)
                                      \big)_\H \big| \, \dt
\\
  & \leq \lim_{T\to\infty}
         \frac{C^s_\QQ}{2T}
         \int_{-T}^T
             \big\| \GRAD e^{-t \AA} W_1 \big\|_{\H^s}
             \big\| e^{-t \AA} W_2 \big\|_\H
             \big\| e^{-t \AA} W_3 \big\|_\H \, \dt
\\
  & = \lim_{T\to\infty}
      \frac{C^s_\QQ}{2T}
      \int_{-T}^T \| \GRAD W_1 \|_{\H^s}
                  \| W_2 \|_\H \| W_3 \|_\H \, \dt
    = C^s_\QQ \| \GRAD W_1 \|_{\H^s}
              \| W_2 \|_\H \| W_3 \|_\H \,.
\end{aligned}
\end{equation*}
This proves the first bound on $\QQBar$ in (\ref{AOp-4}).  The
second bound on $\QQBar$ in (\ref{AOp-4}) follows similarly from
the second bound on $\QQ$ in (\ref{AOp-3}) and the invariances
(\ref{AOp-5}), thereby proving Lemma \ref{AOp-Lem-2}. \qed

\subsection{Averaged Quadratic Convection Operator}

   Until now the entropy has played no role in our analysis of the
weakly nonlinear-dissipative approximation.  It will now play a
central role.  Specifically, the fact system (\ref{HP-1}) has a
strictly convex entropy implies that its flux $F(U)$ has a special
form \cite{Dafermos, Godunov}.  We will use this form to show that
the averaged quadratic convection operator satisfies a cyclic
identity and some new bounds.

\subsubsection{Special Form of the Flux}
Henceforth we will assume that the hyperbolic-parabolic system
(\ref{HP-1}) has a strictly convex entropy that is thrice continuously
differentiable over the convex domain $\UU$.  Define the set
$\VV=\{H_U(U)\,:\,U\in\UU\}$.  One can show that $\VV\subset\RN$ is a
domain and that the mapping $H_U:\UU\to\VV$ is a diffeomorphism.  For
every $V\in\VV$ we then define $H^*(V)\in\R$ and $J^*(V)\in\Rd$ by
\begin{equation}
  \label{AOp-6}
  H^*(V) = V^T U - H(U) \,, \qquad
  J^*(V) = V^T F(U) - J(U) \,,
\end{equation}
where $U\in\UU$ is uniquely determined by $H_U(U)=V$.  Because the
entropy flux $J(U)$ is related to $H(U)$ and $F(U)$ by (\ref{HP-5}),
one can show that the mappings $H^*:\VV\to\R$ and $J^*:\VV\to\Rd$
are continuously differentiable with
\begin{equation}
  \label{AOp-7}
  H^*_V(V) = U \,, \qquad
  J^*_V(V) = F(U) \,.
\end{equation}
It follows that $F(U)$ is given by
\begin{equation}
  \label{AOp-8}
  F(U) = J^*_V\big( H_U(U) \big) \quad
  \hbox{for every $U\in\UU$} \,.
\end{equation}
Because we have assumed that $H$ is thrice continuously
differentiable while $F$ is twice continuously differentiable over
$\UU$, it follows from (\ref{AOp-7}) that both $H^*$ and $J^*$ are
thrice continuously differentiable over $\VV$.  This contrasts with
$J$, which is twice continuously differentiable over $\UU$.

\subsubsection{Cyclic Identity}
The fact that the flux $F(U)$ of system (\ref{HP-1}) has the form
(\ref{AOp-7}) is central to our proof of the following identity.

\begin{Lem}
\label{AOp-Lem-3} {\bf (Cyclic Identity.)} If system (\ref{HP-1})
considered over the spatial domain $\Td$ has a strictly convex
entropy $H$ satisfying (\ref{HP-4}) that is thrice continuously
differentiable over $\UU$ then for every $W_1,W_2,W_3\in
C^\infty(\Td)$ one has the cyclic identity
\begin{equation}
  \label{AOp-9}
  \big( W_1 \,|\, \QQBar(W_2,W_3) \big)_\H
  + \big( W_2 \,|\, \QQBar(W_3,W_1) \big)_\H
  + \big( W_3 \,|\, \QQBar(W_1,W_2) \big)_\H = 0 \,,
\end{equation}
where $\QQBar$ is defined by (\ref{WND-14}).
\end{Lem}

\smallskip
\noindent
{\bf Proof.} We will first establish (\ref{AOp-9}) for the particular
$\QQBar$ associated with the conserved densities.  In that case
$\WW=\UU$.  The extension of (\ref{AOp-9}) to general $\QQBar$ then
follows from the change of variable formula (\ref{WND-18}).

   We will begin by establishing (\ref{AOp-9}) when
$U_1,U_2,U_3\in\oplus_\omega C_\omega^\infty$.  Once this is done,
the extension of (\ref{AOp-9}) to $U_1,U_2,U_3\in C^\infty(\Td)$ then
follows by a density argument that uses the bounds on $\QQBar$ from
Lemma \ref{AOp-Lem-2}.

   Let $U_1,U_2,U_3\in\oplus_\omega C_\omega^\infty$.  From the
definition of $\QQBar$ given by (\ref{WND-14}) we see that
\begin{equation}
  \label{AOp-10}
\begin{aligned}
  \big( U_1 \,|\, \QQBar(U_2,U_3) \big)_\HBbb
   & = \lim_{T\to\infty}
           \frac{1}{2 T}
           \int_{-T}^T
               \big( U_1 \,|\,
                     e^{t\AA} \QQ(e^{-t\AA} U_2, e^{-t\AA} U_3)
                     \big)_\HBbb \, \dt
\\
  & = \lim_{T\to\infty}
          \frac{1}{2 T}
          \int_{-T}^T
              \big( \UTld_1(t) \,|\, \QQ(\UTld_2(t),\UTld_3(t))
                    \big)_\HBbb \, \dt \,.
\end{aligned}
\end{equation}
where $\UTld_1(t)=e^{-t\AA}U_1$, $\UTld_2(t)=e^{-t\AA}U_2$, and
$\UTld_3(t)=e^{-t\AA}U_3$ are quasiperiodic functions of $t$.

   By (\ref{WND-12}) and (\ref{WND-14}) we see that
\begin{equation}
  \nonumber
\begin{aligned}
  \big( \UTld_1 \,|\, \QQ(\UTld_2,\UTld_3) \big)_\HBbb
  & = \int_\Td \UTld_1^T G_o
               \DIV Q_o(\UTld_2,\UTld_3) \, \dx
    = - \int_\Td \GRAD \UTld_1^T
            \DOT G_o Q_o(\UTld_2,\UTld_3) \, \dx \,,
\end{aligned}
\end{equation}
where $G_o=H_{UU}(U_o)$ and
$Q_o(\UTld_2,\UTld_3)=\frac12 F_{UU}(U_o)(\UTld_2\otimes\UTld_3)$.
By differentiating (\ref{AOp-8}) twice and evaluating at $U_o$ we
obtain
\begin{equation}
  \nonumber
\begin{aligned}
  F_U(U_o)
  & = J^*_{VV}(V_o) G_o \,,
\\
  F_{UU}(U_o)\big( \UTld_2 \otimes \UTld_3 \big)
  & = J^*_{VVV}(V_o)\big( G_o \UTld_2 \otimes G_o \UTld_3 \big)
      + J^*_{VV}(V_o) H_{UUU}(U_o)\big( \UTld_2 \otimes \UTld_3 \big) \,.
\end{aligned}
\end{equation}
where $V_o=H_U(U_o)$.  Hence,
\begin{equation}
  \nonumber
  \AA = A_o \DOT \GRAD = J^*_{VV}(V_o) G_o \DOT \GRAD \,,
\end{equation}
and
\begin{equation}
  \nonumber
\begin{aligned}
  \GRAD \UTld_1^T \DOT G_o Q_o(\UTld_2,\UTld_3)
  & = \tfrac12 \GRAD \UTld_1^T \DOT G_o
      J^*_{VVV}(V_o)\big( G_o \UTld_2 \otimes G_o \UTld_3 \big)
\\
  & \quad \,
      + \tfrac12 \GRAD \UTld_1^T \DOT G_o
        J^*_{VV}(V_o) H_{UUU}(U_o)\big( \UTld_2 \otimes \UTld_3 \big)
\\
  & = \tfrac12
      J^*_{VVV}(V_o)\big( G_o \DOT \GRAD \UTld_1
                          \otimes G_o \UTld_2
                          \otimes G_o \UTld_3 \big)
      + \tfrac12
        H_{UUU}(U_o)\big( \del_t \UTld_1
                           \otimes \UTld_2 \otimes \UTld_3 \big) \,.
\end{aligned}
\end{equation}
Summing the above relation with its cyclic permutations gives
\begin{equation}
  \nonumber
\begin{aligned}
  & \GRAD \UTld_1^T \DOT G_o Q_o(\UTld_2,\UTld_3)
    + \GRAD \UTld_2^T \DOT G_o Q_o(\UTld_3,\UTld_1)
    + \GRAD \UTld_3^T \DOT G_o Q_o(\UTld_1,\UTld_2)
\\
  & = \tfrac12
      \DIV \big[ J^*_{VVV}(V_o)\big( G_o \UTld_1 \otimes G_o \UTld_2
                                     \otimes G_o \UTld_3 \big) \big]
      + \tfrac12
        \del_t \big[ H_{UUU}(U_o)\big( \UTld_1 \otimes \UTld_2
                                       \otimes \UTld_3 \big) \big] \,.
\end{aligned}
\end{equation}
Integrating this over $\Td$ we obtain
\begin{equation}
  \nonumber
\begin{aligned}
  & \big( \UTld_1 \,|\, \QQ(\UTld_2,\UTld_3) \big)_\H
    + \big( \UTld_2 \,|\, \QQ(\UTld_3,\UTld_1) \big)_\H
    + \big( \UTld_3 \,|\, \QQ(\UTld_1,\UTld_2) \big)_\H
\\
  & = - \tfrac12
      \frac{\dee}{\dt}
      \int_\Td
          H_{UUU}(U_o)\big( \UTld_1 \otimes \UTld_2
                                    \otimes \UTld_3 \big) \, \dx \,.
\end{aligned}
\end{equation}
The time average of this equation yields
\begin{equation}
  \nonumber
  \big( U_1 \,|\, \QQBar(U_2,U_3) \big)_\H
  + \big( U_2 \,|\, \QQBar(U_3,U_1) \big)_\H
  + \big( U_3 \,|\, \QQBar(U_1,U_2) \big)_\H = 0 \,,
\end{equation}
because the time average of the time derivative of a bounded function
vanishes, thereby proving the lemma. \qed

\smallskip

   The cyclic identity (\ref{AOp-9}) yields the following bound.

\begin{Lem}
\label{AOp-Lem-4}  If system (\ref{HP-1}) considered over the spatial
domain $\Td$ has a strictly convex entropy $H$ satisfying (\ref{HP-4})
that is thrice continuously differentiable over $\UU$ then for every
$s>d/2$ and every $W_1,W_2\in C^\infty(\Td)$ one has the bound
\begin{equation}
  \label{AOp-11}
  \big| \big( W_1 \,|\, \QQBar(W_1,W_2) \big)_\HBbb \big|
  \leq \tfrac12 C^s_\QQ \| W_1 \|_\H^{\,2} \| \GRAD W_2 \|_{\H^s} \,,
\end{equation}
where $\QQBar(W_1,W_2)$ is defined by (\ref{WND-14}) and $C^s_\QQ$ is
the constant appearing in (\ref{AOp-3}) of Lemma \ref{AOp-Lem-1}.
\end{Lem}

\smallskip
\noindent
{\bf Proof.} The cyclic identity (\ref{AOp-9}) and the fact that
$\QQBar(W_1,W_2)=\QQBar(W_2,W_1)$ imply
\begin{equation}
  \nonumber
  2 \big( W_1 \,|\, \QQBar(W_1,W_2) \big)_\HBbb
  + \big( W_2 \,|\, \QQBar(W_1,W_1) \big)_\HBbb = 0 \,.
\end{equation}
It thereby follows from the second bound in (\ref{AOp-4}) that
\begin{equation}
  \nonumber
  \big| \big( W_1 \,|\, \QQBar(W_1,W_2) \big)_\HBbb \big|
  \leq \tfrac12 \big| \big( W_2 \,|\, \QQBar(W_1,W_1) \big)_\HBbb \big|
  \leq \tfrac12 C^s_\QQ \| W_1 \|_\H^{\,2} \| \GRAD W_2 \|_{\H^s} \,,
\end{equation}
where $C^s_\QQ$ is the constant appearing in (\ref{AOp-3}) of
Lemma \ref{AOp-Lem-1}. \qed

%%%%%%%%%%%%%%%%%%%%%%%%%%%%%%%%%%%%%%%%%%%%%%%%%%%%%%%%%%%%%%%%%%%%%%

\subsection{Averaged Dissipation Operator}

   Because the averaged dissipation operator $\DDBar$ is nonpositive
semidefinite, its Friedrichs extension is a self adjoint nonpositive
semidefinite operator with domain $\Dom(\DDBar)$.  For every $V\in\H$
and every $\omega\in\sig(\AA)$ define $V_\omega$ to be the component
of $V$ in $\H_\omega$ --- so that
\begin{equation*}
  V_\omega = E_\omega V \in \H_\omega \,, \qquad
  V = \sum_{\omega\in\sig(\AA)} V_\omega \,.
\end{equation*}
We claim that the nonpositve definite Hermitian form associated with
$\DDBar$ extends to the domain $\Herm(\DDBar)$ given by
\begin{equation}
  \label{AOp-12}
  \Herm(\DDBar)
  = \bigg\{ V \in \H \,:\,
            - \sum_{\omega\in\sig(\AA)}
                  (V_\omega \,|\,\DD V_\omega)_\H
              < \infty \bigg\} \,.
\end{equation}
Recall that $V_\omega\in\Dom(\DD)$ for every $\omega\in\sig(\AA)$
and that $\DD$ is nonpositive over $\Dom(\DD)$, whereby
$(V_\omega \,|\,\DD V_\omega)_\H\leq0$ for every
$\omega\in\sig(\AA)$.  Then for every $V\in\Herm(\DDBar)$ one has
by (\ref{AOp-1}) that
\begin{equation}
  \label{AOp-13}
  \big( V \,|\, \DDBar V \big)_\H
  = \sum_{\omega\in\sig(\AA)}
        (V_\omega \,|\,\DD V_\omega)_\H \ \leq \ 0 \,.
\end{equation}
%The Cauchy-Schwarz inequality then implies that
%\begin{equation*}
%\begin{aligned}
%  \sum_{\omega\in\sig(\AA)} (V_\omega \,|\,\DD V_\omega)_\H
%  = \big( V \,|\, \DDBar V \big)_\H
%  & \leq \| V \|_\H \big\| \DDBar V \big\|_\H
%    = \| V \|_\H
%      \Bigg( \sum_{\omega\in\sig(\AA)}
%                 \| E_\omega \DD V_\omega \|_{\H}^{\,2}
%             \Bigg)^\frac12 \,.
%\end{aligned}
%\end{equation*}
%It follows that $\Dom(\DDBar)\subset\Herm(\DDBar)\subset\H$.

   We now characterize when the Hermitian form associated with
$\DDBar$ is nondegenerate.

\begin{Lem}
  \label{AOp-Lem-5}
\begin{equation}
  \label{AOp-14}
  \big( V \,|\, \DDBar V \big)_\H < 0 \quad
  \hbox{for every nonconstant $V\in\Herm(\DD)$} \,,
\end{equation}
if and only if the Kawashima condition (\ref{WND-7}) holds.
\end{Lem}

\smallskip
\noindent
{\bf Remark.}  This kind of characterization was first proved by
Kawashima \cite{Ka,KaSh}.

\smallskip
\noindent
{\bf Proof.} First we show that (\ref{AOp-14}) implies the Kawashima
condition (\ref{WND-7}) holds.  Suppose not.  Then there exists a
nonconstant $V\in\H$ such that $V\in\H_\omega$ for some
$\omega\in\sig(\AA)$ and $V\in\Null(\DD)$.  But then
$V=V_\omega\in\Herm(\DDBar)$ with
\begin{equation*}
  \big( V \,|\, \DDBar V \big)_\H
  = (V_\omega \,|\,\DD V_\omega)_\H
  = (V \,|\,\DD V)_\H = 0 \,,
\end{equation*}
which contradicts (\ref{AOp-14}).  Therefore (\ref{AOp-14}) implies
the Kawashima condition (\ref{WND-7}) holds.

   Next we show that the Kawashima condition (\ref{WND-7}) implies
(\ref{AOp-14}).  Suppose that (\ref{AOp-14}) is false.  Then there
exists a nonconstant $V\in\Herm(\DDBar)$ such that
\begin{equation*}
  0 = \big( V \,|\, \DDBar V \big)_\H
    = \sum_{\omega\in\sig(\AA)} (V_\omega \,|\,\DD V_\omega)_\H \,.
\end{equation*}
However this is equivalent to
\begin{equation*}
  (V_\omega \,|\,\DD V_\omega)_\H = 0 \quad
  \hbox{for every $\omega\in\sig(\AA)$} \,,
\end{equation*}
which is equivalent to
\begin{equation*}
  V_\omega \in \Null(\DD) \quad
  \hbox{for every $\omega\in\sig(\AA)$} \,,
\end{equation*}
But $V\neq0$ implies that $V_\omega\neq0$ for some
$\omega\in\sig(\AA)$.  But then for this $\omega$ we have $V_\omega$
is nonconstant, $V_\omega\in\H_\omega$, and $V_\omega\in\Null(\DD)$,
which contradicts the Kawashima condition (\ref{WND-7}).  Therefore
the Kawashima condition (\ref{WND-7}) implies (\ref{AOp-14}), and the
proof of Lemma \ref{AOp-Lem-5} is complete.
\qed

%%%%%%%%%%%%%%%%%%%%%%%%%%%%%%%%%%%%%%%%%%%%%%%%%%%%%%%%%%%%%%%%%%%%%%
%%  Section 5: Strictly Dissipative Approximations                  %%
%%%%%%%%%%%%%%%%%%%%%%%%%%%%%%%%%%%%%%%%%%%%%%%%%%%%%%%%%%%%%%%%%%%%%%

\section{Strictly Dissipative Approximations}
\label{SDA}

   The global existence theory presented in the next section will
require that $\DDBar$ be strictly elliptic.  When this is the
case we say that the WND approximation is {\em strictly dissipative}.
For many physical systems the WND approximation has this property.
Indeed, the averaged dissipation operator has played an important
role in the study of large time behavior of solutions to
hyperbolic-parabolic systems of conservations laws and discrete
Boltzmann equations.  Kawashima \cite{Kawashima, Kawashima-BE}
showed that in one space dimension the averaged dissipation
operator was strictly dissipative whenever the Kawashima condition
(\ref{WND-7}) holds.  In that case he showed that solutions of the
original system are well approximated by solutions to an
``effective artificial viscosity'' system constructed using the
averaged dissipation operator.  Motivated by this idea, Hoff and
Zumbrun \cite{HZ1, HZ2} studied multi-dimensional diffusion waves
for the barotropic Navier-Stokes system through an artificial
viscosity system constructed with the averaged dissipation operator.
Recently, Bianchini-Hanouzet-Natalini \cite{BHN} used the same idea
to study the large-time behavior of smooth solutions for partially
dissipative hyperbolic systems with a convex entropy.  In each case
the averaged dissipation operator was shown to be strictly
dissipative through a detailed spectral analysis.

   To our knowledge there is no proof that the Kawashima condition
(\ref{WND-7}) implies that the averaged dissipation operator is
strictly dissipative in multidimensional settings.  Here we give
a stronger criterion that does the job in our spatially periodic
setting.  In the Fourier representation we have
\begin{equation*}
  \widehat{\AA V}(\xi) = i A_o(\xi) \VHat(\xi) \,, \qquad
  \widehat{\DD V}(\xi) = - B_o(\xi) \VHat(\xi) \,,
\end{equation*}
where $A_o(\xi)$ and $B_o(\xi)$ are the families of $G_o$-symmetric
matrices in $\RNN$ defined for every $\xi\in\Rd$ by
\begin{equation*}
  A_o(\xi) = \xi \DOT A_o \,, \qquad
  B_o(\xi) = \xi^{\otimes2} \DDOT B_o \,.
\end{equation*}
One then has
\begin{equation*}
  \widehat{e^{-t\AA} V}(\xi) = e^{-itA_o(\xi)} \VHat(\xi) \,,
\end{equation*}
whereby
\begin{equation*}
  \widehat{\DDBar V}(\xi) = \DDBarHat(\xi) \VHat(\xi) \,, \quad
  \DDBarHat(\xi)
  = \lim_{T\to\infty}
        \frac{1}{2 T}
        \int_{-T}^T e^{i t A_o(\xi)} B_o(\xi)
                    e^{-i t A_o(\xi)} \, \dt \,.
\end{equation*}

   Our Kawashima-type criterion for strict dissipativity is given
by the following.

\begin{Thm}
  \label{SDA-Thm}
Suppose that for some $\alpha>0$ there exists $\beta>0$ such that
\begin{equation}
  \label{SDA-1}
  G_o B_o(\xiHat)
  + \frac{1}{\alpha^2} \,
    A_o(\xiHat)^T G_o B_o(\xiHat) A_o(\xiHat)
  \geq \beta \, G_o \quad
  \hbox{for every $\xiHat\in\Sd$} \,.
\end{equation}
Then there exists $\delta>0$ such that
\begin{equation}
  \label{SDA-2}
  - \big( V \,|\, \DDBar V \big)_\H
  \geq \delta \, \| \GRAD V \|_\H^2 \quad
  \hbox{for every $V\in C^2(\Td)$} \,.
\end{equation}
\end{Thm}

\smallskip
\noindent
{\bf Remark.} The Kawashima condition is satisfied whenever
(\ref{SDA-1}) holds.

\smallskip

   Our proof of Theorem \ref{SDA-Thm} requires the following lemma.

\begin{Lem}
  \label{SDA-Lem}
Let $G\in\CNN$ be Hermitian positive definite (i.e. $G=G^*>0$).  Let
$A\in\CNN$ be $G$-Hermitian (i.e. $GA=A^*G$) and $B\in\CNN$ be
$G$-Hermitian nonnegative definite (i.e. $GB=B^*G\geq0$).  Let $C_A$
and $C_B$ be constants such that $\|A\|_G\leq C_A$ and
$\|B\|_G\leq C_B$.  Suppose that for some $\alpha>0$ there exists
$\beta>0$ such that
\begin{equation}
  \label{SDA-3}
  G B + \frac{1}{\alpha^2} \, A^* G B A \geq \beta \, G \,.
\end{equation}
Then there exists $\delta>0$ depending only on $\alpha$, $\beta$,
$C_A$, and $C_B$, such that
\begin{equation}
  \label{SDA-4}
  \lim_{T\to\infty}
      \frac{1}{2T} \int_{-T}^T e^{i t A^*} G B e^{- i t A} \dt
  \geq \delta \, G  \,.
\end{equation}
\end{Lem}

\smallskip
\noindent
{\bf Proof.}  Let $v\in\CN$ be nonzero and define
\begin{equation*}
  f(t) = v^* e^{i t A^*} G B e^{- i t A} v \,.
\end{equation*}
This is a smooth, nonnegative, quasiperiodic, real-valued function
that will vanish whenever $e^{-itA}v$ is a null vector of $B$.  In
order to prove (\ref{SDA-4}) we must obtain a lower bound for its
average.  The idea of the proof is to show that $f(t)$ cannot be
too small for long.

   The first two derivatives of $f(t)$ are
\begin{equation}
  \label{SDA-5}
\begin{aligned}
  \dot{f}(t)
  & = i v^* e^{i t A^*} (A^* G B - G B A) e^{- i t A} v \,,
\\
  \ddot{f}(t)
  & = v^* e^{i t A^*} (2 A^* G B A - A^{*2} G B
                                   - G B A^2) e^{- i t A} v \,.
\end{aligned}
\end{equation}
Roughly speaking, we will show that $\ddot{f}(t)$ is dominated by
its first term when $f(t)$ is small.  For every $\eta>0$ one has
\begin{equation*}
\begin{aligned}
  \big| v^* e^{i t A^*} (A^{*2} G B + G B A^2) e^{- i t A} v \big|
  & \leq \frac{1}{\eta^2} \, v^* e^{i t A^*} G B e^{- i t A} v
         + \eta^2 v^* e^{i t A^*} A^{*2} G B A^2 e^{- i t A} v \,,
\end{aligned}
\end{equation*}
The first term on the right-hand side above is just $f(t)/\eta^2$
while the second can be bounded by $\eta^2 C_A^{\,4}C_B^{}v^*Gv$.
If we set $\eta^2=\alpha^2\beta/(2 C_A^{\,4}C_B^{})$ then we obtain
the bound
\begin{equation}
  \label{SDA-6}
  \big| v^* e^{i t A^*} (A^{*2} G B + G B A^2) e^{- i t A} v \big|
  \leq \frac{2 C_A^{\,4} C_B^{}}{\alpha^2 \beta} \, f(t)
       + \frac{\alpha^2\beta}{2} \, v^* G v \,.
\end{equation}

   Let $\Omega^\eps=\{t\in\R\,:\,f(t)<\eps\,v^*Gv\}$, where $\eps>0$
satisfies
\begin{equation}
  \label{SDA-7}
  \eps \leq \frac{\alpha^4 \beta}{\alpha^4 \beta + C_A^{\,4} C_B^{}} \,
            \frac{\beta}{4} \,.
\end{equation}
For every $t\in\Omega^\eps$ we obtain from (\ref{SDA-3}),
(\ref{SDA-5}), (\ref{SDA-6}), and (\ref{SDA-7}) the lower bound
\begin{equation}
  \label{SDA-8}
\begin{aligned}
  \ddot{f}(t)
  & \geq v^* e^{i t A^*} (2 \alpha^2 G B + 2 A^* G B A) e^{- i t A} v
         - \bigg( 2 \alpha^2 + \frac{2 C_A^{\,4} C_B^{}}
                                    {\alpha^2 \beta} \bigg) f(t)
         - \frac{\alpha^2 \beta}{2} \, v^* G v
\\
  & \geq \bigg( 2 \alpha^2 \beta
                - \bigg( 2 \alpha^2
                         + \frac{2 C_A^{\,4} C_B^{}}
                                {\alpha^2 \beta} \bigg) \eps
                - \frac{\alpha^2\beta}{2} \bigg) v^* G v
\\
  & \geq \alpha^2 \beta \, v^* G v \,.
\end{aligned}
\end{equation}
For every $t\in\Omega^\eps$ we obtain from (\ref{SDA-5}),
(\ref{SDA-6}), and (\ref{SDA-7}) the upper bound
\begin{equation}
  \label{SDA-9}
\begin{aligned}
  \ddot{f}(t)
  & \leq v^* e^{i t A^*} (2 A^* G B A) e^{- i t A} v
         + \frac{2 C_A^{\,4} C_B^{}}{\alpha^2 \beta} \, f(t)
         + \frac{\alpha^2 \beta}{2} \, v^* G v
\\
  & \leq \bigg( 2 C_A^{\,2} C_B^{}
                + \frac{2 C_A^{\,4} C_B^{}}{\alpha^2 \beta} \, \eps
                + \frac{\alpha^2 \beta}{2} \bigg) v^* G v
\\
  & \leq \big( 2 C_A^{\,2} C_B^{} + \alpha^2 \beta \big) v^* G v \,.
\end{aligned}
\end{equation}

   Because $f$ is continuous the set $\Omega^\eps$ is open and is
therefore a countable union of disjoint open intervals:
\begin{equation*}
  \Omega^\eps = \bigcup_{k\in\N} (a_k^\eps,b_k^\eps) \,.
\end{equation*}
Let $(a,b)$ be any one of these intervals.  Because $\ddot{f}(t)$
satisfies the lower bound (\ref{SDA-8}) while $f(t)<\eps\,v^*Gv$
for every $t\in(a,b)$, it is clear that the interval $(a,b)$ must
be bounded.

   We will begin by bounding $b-a$ above and below.  Because $f$
is continuous over the bounded interval $[a,b]$, and because $a$
and $b$ are not in $\Omega^\eps$, it follows that
$f(a)=f(b)=\eps\,v^*Gv$, and that $f$ takes its minimum at a
point $t_o\in(a,b)$, at which $\dot{f}(t_o)=0$.  Then
\begin{equation}
  \label{SDA-10}
  f(t) = f(t_o) + \int_{t_o}^t \dot{f}(t_1) \, \dt_1
       = f(t_o) + \int_{t_o}^t \int_{t_o}^{t_1}
                      \ddot{f}(t_2) \, \dt_2 \, \dt_1 \,.
\end{equation}
We claim that $b-a$ satisfies the bounds
\begin{equation}
  \label{SDA-11}
  \frac{\eps \, v^* G v - f(t_o)}
       {2 C_A^{\,2} C_B^{} + \alpha^2 \beta}
  \leq \tfrac18 (b - a)^2 v^* G v
  \leq \frac{\eps \, v^* G v}{\alpha^2 \beta} \,.
\end{equation}

   The upper bound of (\ref{SDA-11}) is obtained from (\ref{SDA-10})
by using the lower bound (\ref{SDA-8}) for $\ddot{f}(t)$.  For every
$t\in(a,b)$ we have
\begin{equation*}
  f(t)
  \geq f(t_o) + \tfrac12 (t - t_o)^2 \alpha^2 \beta \, v^* G v \,.
\end{equation*}
Evaluating this at $t=a$ and $t=b$ yields
\begin{equation*}
\begin{aligned}
  \eps \, v^* G v
  & \geq f(t_o) + \tfrac12 (t_o - a)^2 \alpha^2 \beta \, v^* G v \,,
\\
  \eps \, v^* G v
  & \geq f(t_o) + \tfrac12 (b - t_o)^2 \alpha^2 \beta \, v^* G v \,.
\end{aligned}
\end{equation*}
Because $f(t_o)\geq0$ while
$\max\{(t_o - a)^2,(b - t_o)^2\}\geq\tfrac14(b-a)^2$, we obtain
\begin{equation*}
  \eps \, v^* G v
  \geq \tfrac18 (b - a)^2 \alpha^2 \beta \, v^* G v \,,
\end{equation*}
which yields the upper bound of (\ref{SDA-11}).

   The lower bound of (\ref{SDA-11}) is obtained from (\ref{SDA-10})
by using the upper bound (\ref{SDA-9}) for $\ddot{f}(t)$.  For every
$t\in(a,b)$ we have
\begin{equation*}
  f(t)
  \leq f(t_o) + \tfrac12 (t - t_o)^2
                (2 C_A^{\,2} C_B^{} + \alpha^2 \beta) \, v^* G v \,,
\end{equation*}
Evaluating this at $t=a$ and $t=b$ yields
\begin{equation*}
\begin{aligned}
  \eps \, v^* G v
  & \leq f(t_o) + \tfrac12 (t_o - a)^2
                  (2 C_A^{\,2} C_B^{} + \alpha^2 \beta) \, v^* G v \,,
\\
  \eps \, v^* G v
  & \leq f(t_o) + \tfrac12 (b - t_o)^2
                  (2 C_A^{\,2} C_B^{} + \alpha^2 \beta) \, v^* G v \,.
\end{aligned}
\end{equation*}
Because $\min\{(t_o - a)^2,(b - t_o)^2\}\leq\tfrac14(b-a)^2$, we obtain
\begin{equation*}
  \eps \, v^* G v
  \leq f(t_o) + \tfrac18 (b - a)^2
                (2 C_A^{\,2} C_B^{} + \alpha^2 \beta) \, v^* G v \,,
\end{equation*}
which yields the lower bound of (\ref{SDA-11}).

   Next, we bound the average of $f(t)$ over $(a,b)$ from below.
By again using the lower bound (\ref{SDA-8}) for $\ddot{f}(t)$ in
(\ref{SDA-10}) we obtain
\begin{equation*}
\begin{aligned}
  \frac{1}{b - a} \int_a^b f(t) \, \dt
  & \geq f(t_o)
         + \frac{(b - t_o)^2 + (b - t_o)(t_o - a) + (t_o - a)^2}{6}
           \, \alpha^2 \beta \, v^* G v
\\
  & \geq f(t_o) + \tfrac18 (b - a)^2 \alpha^2 \beta \, v^* G v \,.
\end{aligned}
\end{equation*}
The lower bound  of (\ref{SDA-11}) then implies
\begin{equation}
  \label{SDA-12}
\begin{aligned}
  \frac{1}{b - a} \int_a^b f(t) \, \dt
  & \geq f(t_o) + \alpha^2 \beta \,
                  \frac{\eps \, v^* G v - f(t_o)}
                       {2 C_A^{\,2} C_B^{} + \alpha^2 \beta}
    \geq \frac{\alpha^2 \beta \, \eps}
              {2 C_A^{\,2} C_B^{} + \alpha^2 \beta} \, v^* G v \,.
\end{aligned}
\end{equation}

   We are now ready to prove (\ref{SDA-4}) with $\delta$ given by
\begin{equation}
  \label{SDA-13}
   \delta = \frac{\alpha^2 \beta \, \eps}
                 {2 C_A^{\,2} C_B^{} + \alpha^2 \beta} \,.
\end{equation}
For every $T>0$ let
$K_T^\eps=\{k\in\N\,:\,(a_k^\eps,b_k^\eps)\subset(-T,T)\}$
and
\begin{equation*}
  \Omega_T^\eps = \bigcup_{k\in K^\eps_T} (a_k^\eps,b_k^\eps) \,.
\end{equation*}
Then using the fact that $f(t)\geq\eps$ over $[-T,T]-\Omega^\eps$,
the lower bound (\ref{SDA-12}), and the fact that $\delta$ given by
(\ref{SDA-13}) satisfies $\delta<\eps$,  we find
\begin{equation*}
\begin{aligned}
  \frac{1}{2 T} \int_{-T}^T f(t) \, \dt
  & \geq \frac{1}{2 T}
         \bigg[ \int_{[-T,T]-\Omega^\eps} f(t) \, \dt
                + \sum_{k\in K_T^\eps} \int_{a_k}^{b_k} f(t) \, \dt
                \bigg]
\\
  & \geq \frac{1}{2 T}
         \bigg[ {\mathrm{meas}}\big([-T,T]-\Omega^\eps\big) \, \eps
                + \sum_{k\in K_T^\eps} (b_k^\eps - a_k^\eps) \, \delta
                \bigg] v^* G v
\\
  & \geq \bigg[ 1 - \frac{{\mathrm{meas}}\big( [-T,T] \cap \Omega^\eps
                                              - \Omega_T^\eps \big)}
                         {2 T} \bigg] \delta \, v^* G v \,.
\end{aligned}
\end{equation*}
The set $[-T,T]\cap\Omega^\eps-\Omega_T^\eps$ is contained in the union
of the (at most two) disjoint intervals $(a^\eps_k,b^\eps_k)$ that
contain $-T$ and $T$.  So its measure is bounded above by twice the
upper bound for $b-a$ given by (\ref{SDA-11}) --- namely, by
\begin{equation*}
  {\mathrm{meas}}\big( [-T,T] \cap \Omega^\eps - \Omega_T^\eps \big)
  \leq 4 \bigg( \frac{2 \eps}{\alpha^2 \beta} \bigg)^\frac12 \,.
\end{equation*}
Hence, we obtain
\begin{equation*}
  \frac{1}{2 T} \int_{-T}^T f(t) \, \dt
  \geq \bigg[ 1 - \frac{2}{T}
                  \bigg( \frac{2 \eps}{\alpha^2 \beta} \bigg)^\frac12
              \bigg]  \delta \, v^* G v \,.
\end{equation*}
Letting $T\to\infty$ in this inequality yields (\ref{SDA-4}).  The
limit of the left-hand side exists because $f(t)$ is quasiperiodic.
\qed

\medskip

   We are now ready to prove Theorem \ref{SDA-Thm} with the aid of
Lemma \ref{SDA-Lem}.

\smallskip
\noindent
{\bf Proof of Theorem \ref{SDA-Thm}.}  First observe that for every
nonzero $\xi\in\Rd$ we set $\xiHat=\xi/|\xi|$ and have
\begin{equation*}
\begin{aligned}
  - \DDBarHat(\xi)
  & = \lim_{T\to\infty}
          \frac{1}{2 T}
          \int_{-T}^T e^{i t A_o(\xi)} B_o(\xi)
                      e^{-i t A_o(\xi)} \, \dt
\\
  & = |\xi|^2
      \lim_{T\to\infty}
          \frac{1}{2 T}
          \int_{-T}^T e^{i t |\xi| A_o(\xiHat)} B_o(\xiHat)
                      e^{-i t |\xi| A_o(\xiHat)} \, \dt
\\
  & = |\xi|^2
      \lim_{T\to\infty}
          \frac{1}{2 T}
          \int_{-T}^T e^{i t A_o(\xiHat)} B_o(\xiHat)
                      e^{-i t A_o(\xiHat)} \, \dt \,.
\end{aligned}
\end{equation*}
Then by the fact $A_o(\xiHat)$ is $G_o$-symmetric we have
\begin{equation}
  \label{SDA-14}
\begin{aligned}
  - G_o \DDBarHat(\xi)
  & = |\xi|^2
      \lim_{T\to\infty}
          \frac{1}{2 T}
          \int_{-T}^T G_o e^{i t A_o(\xiHat)} B_o(\xiHat)
                          e^{-i t A_o(\xiHat)} \, \dt
\\
  & = |\xi|^2
      \lim_{T\to\infty}
          \frac{1}{2 T}
          \int_{-T}^T e^{i t A_o(\xiHat)^T} G_o B_o(\xiHat)
                      e^{-i t A_o(\xiHat)} \, \dt \,.
\end{aligned}
\end{equation}

   We now apply Lemma \ref{SDA-Lem} with $G=G_o$, $A=A_o(\xiHat)$,
$B=B_o(\xiHat)$,
\begin{equation*}
  C_A = \max\big\{ \| A_o(\xiHat) \|_{G_o} \,:\,
                   \xiHat \in \Sd \big\} \,, \qquad
  C_B = \max\big\{ \| B_o(\xiHat) \|_{G_o} \,:\,
                   \xiHat \in \Sd \big\} \,.
\end{equation*}
We find that for every $\xiHat\in\Sd$ we have the lower bound
\begin{equation*}
  \lim_{T\to\infty}
      \frac{1}{2 T}
      \int_{-T}^T e^{i t A_o(\xiHat)^T} G_o B_o(\xiHat)
                  e^{-i t A_o(\xiHat)} \, \dt
  \geq \delta \, G_o \,,
\end{equation*}
where $\delta$ is given by (\ref{SDA-13}).  Because $\delta$ only depends
on $\alpha$, $\beta$, $C_A$, and $C_B$, it is independent of $\xiHat$.
Combining this lower bound with (\ref{SDA-14}) yields
\begin{equation*}
  - G_o \DDBarHat(\xi) \geq |\xi|^2 \delta \, G_o \quad
  \hbox{for every $\xi\in\Rd$} \,.
\end{equation*}

   Then for every $V\in C^2(\Td)$ the Plancherel identity implies
\begin{equation*}
\begin{aligned}
  - \big( V \,|\, \DDBar V \big)_\H
  & = - \sum_{\xi\in\Ld^*}
            \VHat(\xi)^* G_o \DDBarHat(\xi) \VHat(\xi)
    \geq \delta \, \sum_{\xi\in\Ld^*}
                       |\xi|^2 \VHat(\xi)^* G_o \VHat(\xi)
    = \delta \, \| \GRAD V \|_\H^2 \,.
\end{aligned}
\end{equation*}
But this is (\ref{SDA-2}), thereby proving the theorem.  \qed

%%%%%%%%%%%%%%%%%%%%%%%%%%%%%%%%%%%%%%%%%%%%%%%%%%%%%%%%%%%%%%%%%%%%%%
%%  Section 6: Global Weak Solutions                                %%
%%%%%%%%%%%%%%%%%%%%%%%%%%%%%%%%%%%%%%%%%%%%%%%%%%%%%%%%%%%%%%%%%%%%%%

\section{Global Weak Solutions}
\label{GWS}

   The weakly nonlinear-dissipative approximation is
\begin{equation}
  \label{GWS-1}
  \del_t W + \AA W + \QQBar(W,W)
  = \DDBar W \,.
\end{equation}
Following the Leray theory for the incompressible Navier-Stokes
system, we will show that if $\DDBar$ is strictly dissipative
then (\ref{GWS-1}) has global weak solutions for all initial data
$W^\init\in\HBbb$.  This result includes the Leray theory, so it
cannot be improved easily.

   The key to obtaining global solutions in the Leray theory for the
incompressible Navier-Stokes system is a so-called energy estimate.
This designation is a bit misleading because, as we shall see, the
estimate is better understood as an entropy estimate.

\subsection{Notion of Weak Solution}

   We call $W \in C([0,\infty);\wH)\cap L^2_\loc(\dt;\VBbb)$ a
Leray-type weak solution of weakly nonlinear-dissipative
approximation \eqref{GWS-1}, if $W$ satisfies the following weak
form of \eqref{GWS-1}:
\begin{equation}
  \label{GWS-2}
\begin{aligned}
  0 & = \big( V \,|\, W(t_2) \big)_{\HBbb}
        - \big( V \,|\, W(t_1) \big)_{\HBbb}
        - \int_{t_1}^{t_2}
          \big( \AA V \,|\, W(t) \big)_{\HBbb} \, \dt
\\
  & \quad \,
      - \int_{t_1}^{t_2} \!\! \int_\Td
            \GRAD V
            \DOT G_o \QBar(W(t),W(t)) \, \dx \, \dt
      + \int_{t_1}^{t_2} \!\! \int_\Td
            \GRAD V
            \DOT G_o \BBar \DOT \GRAD W(t) \, \dx \, \dt \,,
\end{aligned}
\end{equation}
for every function $V \in \V$.  These solutions satisfy the
entropy inequality
\begin{equation}
  \label{GWS-3}
  \tfrac12 \| W(t) \|_{\HBbb}^{\,2}
  - \int_0^t (W \,|\, \DDBar W)_{\HBbb} \, \dt
  \leq \tfrac12 \| W^\init \|_{\HBbb}^{\,2}\,.
\end{equation}
Of course, for every sufficiently nice $W$  one has the
identities
\begin{equation}
  \label{Identity}
  \big( W \,|\, \AA W \big)_{\HBbb} = 0 \,, \qquad
  \big( W \,|\, \QQBar(W, W) \big)_{\HBbb} = 0 \,.
\end{equation}

\subsection{Existence Theorem}

    The main theorem of this paper, the global existence of
Leray-type weak solution to the weakly nonlinear-dissipative
approximation \eqref{GWS-1} is as follows.

\begin{Thm}
Let \eqref{GWS-1} be the weakly nonlinear-dissipative approximation
about a constant state $U_o$ of the hyperbolic-parabolic system
(\ref{HP-1}) with a strictly convex entropy (\ref{HP-4}).  Suppose
that the diffusive operator $\DBar$ in \eqref{GWS-1} is strictly
dissipative, i.e. that there exists $\delta>0$ such that
\begin{equation}
  \label{dissipate-D}
  - \big( V \,|\, \DDBar V \big)_\H
  \geq \delta \, \| \GRAD V \|_\H^2 \quad
  \hbox{for every $V\in C^2(\Td)$}\,.
\end{equation}
Then, for every $\WTld^\init \in \V$, there exists a solution
$\WTld\in C([0,\infty);\wH)\cap L^2_\loc(\dt;\V)$ satisfying
\eqref{GWS-2} and \eqref{GWS-3}.
\end{Thm}

\smallskip
\noindent
{\bf Remark.} If condition \eqref{SDA-1} is satisfied then $\DDBar$
will satisfy \eqref{dissipate-D} and the above theorem insures the
existence of at least one global weak solution.

\smallskip
\noindent
{\bf Remark.}  In Section \ref{WCNS}, we apply our theory to the
Navier-Stokes system of gas dynamics.  The resulting averaged system
includes the incompressible Navier-Stokes system as a subsystem.
The question of uniqueness for system \eqref{GWS-1} is thereby at
least as hard as that of uniqueness for weak solutions of the
incompressible Navier-Stokes system.

\smallskip
\noindent {\bf Proof.}  The strategy for our proof was introduced by
Leray in the context of the incompressible Navier-Stokes system
\cite{Le34}, see also \cite{La, Te}. It is now classical compactness
argument that has since been used to prove existence of global weak
solutions for other equations \cite{DGL}. It proceeds in four steps.
We begin by constructing a sequence of approximate solutions.  We
then show that this sequences is relatively compact, first in some
weak topologies and then in a strong topology. Finally, we show that
limit points of this sequence satisfy (\ref{GWS-2}) and are thereby
weak solutions of (\ref{GWS-1}).  This strategy strikes a balance
between the fact that compactness is easier to establish for weaker
topologies and the fact that passing to the limit in nonlinear terms
requires convergence in a strong topology.

\subsubsection{Step 1: Constructing Approximate Solutions}

   One can construct a sequence of approximation solutions $W_n$
by any method that yields a consistent weak formulation and an
energy relation.  Here we do this with the Galerkin method.

   Let $\{\H_n\}_{n=1}^\infty$ be a sequence of subspaces of $\H$
such that each $\H_n$ lies within $C^\infty(\Td)$, has dimension
$n$, and satisfies $\H_n\subset\H_{n+1}$.  Assume moreover that
this sequence is complete.  Let $\PP_n$ denote the orthogonal
projection from $\H$ onto $\H_n$.  Completeness implies that for
every $V\in\H$ one has $\PP_nV\to V$ as $n\to\infty$.  The
Galerkin approximation of dimension $n$ is the system
\begin{equation}
  \label{Galerkin}
  \del_t W_n + \PP_n \AA W_n + \PP_n \QQBar(W_n, W_n)
  = \PP_n \DDBar W_n\,,
\end{equation}
where $W_n$ takes values in $\H_n$.  This is a system of $n$ ODEs.
Its nonlinearities are quadratic, hence locally Lipschitz.  The
Picard existence theorem insures that system \eqref{Galerkin} has
local solutions. Taking inner product with $W_n$, and applying the
identities in \eqref{Identity}, we obtain the energy identity:
\begin{equation}
  \label{Galerkin-energy}
  \tfrac12 \| W_n(t) \|_\H^{\,2}
  - \int_0^t (W_n \,|\, \DDBar W_n)_\H \, \dt'
  = \tfrac12 \| W_n^\init \|_\H^{\,2}\,,
\end{equation}
for every $t>0$.  This energy identity immediately implies a
global $L^2$ bound on the approximate solutions $W_n$, which
thereby exists for all time.

\subsubsection{Step 2: Establish Weak Compactness}

   We claim that the approximate solutions $W_n$ are relatively
compact in $C([0,\infty);\wH)\cap\mbox{w-}L^2_\loc(\dt;\wV)$.

   First, from the energy identity \eqref{Galerkin-energy}, we that
$W_n(t)$ is uniformly bounded in $\H$ thus relatively compact in
$\wH$ for every $t>0$.  Next, from the dissipation property of
$\DDBar$, see \eqref{dissipate-D}, we have
\begin{equation}
  \delta \int_0^t \| \GRAD W_n\|^2_\H \, \dt'
  \leq - \int_0^t \big( W_n\, |\, \DDBar W_n\big)_H \,\dt'
  \leq C \,.
 \end{equation}
Thus $W_n$ is relatively compact in
$\mbox{w-}L^2_{loc}([0,\infty);\wV)$.  We need only to verify that
$W_n$ is equicontinuous in $C([0,\infty);\wH)$ thus by Arzela-Ascoli
theorem the weak compactness is established.  The equicontinuity can
be derived from the weak form of the Galerkin system
\eqref{Galerkin}.
\begin{equation}
   \label{Galerkin-weakform}
\begin{aligned}
  \big( \VTld \,|\, W_n(t_2) - W_n(t_1) \big)_\H
  & =  \int^{t_2}_{t_1} \big( \VTld \,|\, \AA W_n(t) \big)_\H \, \dt
\\
  & \quad \,
      - \int^{t_2}_{t_1}
           \big( \VTld \,|\, \QQBar( W_n(t), W_n(t))\big)_\H\, \dt
      + \int^{t_2}_{t_1} \big( \VTld \,|\, \DDBar W_n(t)\big)_\H\, \dt \,,
\end{aligned}
\end{equation}
for every function $\VTld\in\V$.  We first prove the
equicontinuity for test function $\VTld \in \V\cap C^1(\T^d)$, which
is followed from the first and third bounds in the Lemma
\ref{AOp-Lem-2}.  Then we extend the class of test functions to $\V$
by standard density argument, thereby finishing the proof of Step 2.

\subsubsection{Step 3: Establish Strong Compactness}

   We claim that $W_n$ is relatively compact in strong topology of
$L^2_{loc}([0,\infty);\H)$.  It is a direct consequence of the
weak compactness result in Step 2 and the fact that the injection
\begin{equation}
  C([0,\infty);\wH)\cap\mbox{w-}L^2_\loc(\dt;\wV)
 \To  L^2_{loc}([0,\infty);\H)
\end{equation}
is continuous.

\subsubsection{Step 4: Pass to the Limit}

   Step 2 ensures that there is a subsequence of $W_n$, which we
also refer to as $W_n$, converges in
$C([0,\infty);\wH)\cap\mbox{w-}L^2_\loc(\dt;\wV)$ to a limit
$ W\in C([0,\infty);\wH)\cap L^2_\loc(\dt;\V)$.  Step 3 ensures the
convergence of $W_n$ to $W$ in $L^2_{loc}([0, \infty); \H)$. All
that remains is to show that the limit $W$ satisfies the weak
form \eqref{GWS-2} as well as the energy inequality \eqref{GWS-3}.
Toward this end we check convergence of each term in the respective
regularized versions, \eqref{Galerkin-weakform} and
\eqref{Galerkin-energy}, respectively. Again, we first consider the
test function $V$ in the class $\V\cap C^1(\T^d)$ then use
density argument later. First
\begin{equation}
  \nonumber
  \big( V \,|\, W_n(t_2) - W_n(t_1) \big)_\H
  \to \big( V \,|\, W(t_2) - W(t_1) \big)_\H
  \quad \mbox{as} \quad\! n \to \infty \,,
\end{equation}
because of the relative compactness of $W_n$ in $C([0,\infty);\wH)$.
The convergence of the first term on the righthand side of
\eqref{Galerkin-weakform} is trivial. Note that
\begin{equation}
  \nonumber
  \int^{t_2}_{t_1}
      \big( V \,|\, \QQBar( W_n(t), W_n(t)) \big)_\H \, \dt
  = - \int_{t_1}^{t_2} \!\! \int_\Td
          \GRAD V
          \DOT G_o \QBar(W_n(t),W_n(t)) \, \dx \, \dt \,,
\end{equation}
and $\QBar(W_n(t),W_n(t))$ is quadratic in $W_n$. Thus the strong
compactness of $W_n$ in $L^2_{loc}([0,\infty);\H)$ in Step 3 ensures
the convergence of above term. We also note that
\begin{equation}
  \nonumber
  \int^{t_2}_{t_1} \big( V \,|\, \DDBar W_n(t)\big)_\H\, \dt
  = \int_{t_1}^{t_2} \!\! \int_\Td
        \GRAD V \DOT G_o \BBar \DOT \GRAD W_n(t) \, \dx \, \dt \,,
\end{equation}
The convergence of above term is straightforward.  Thus, we show that
the limit $W$ satisfy the weak form \eqref{GWS-2}, thus is a
weak solution to WND approximation \eqref{GWS-1}.

   Now, to recover the energy inequality \eqref{GWS-3} from
\eqref{Galerkin-energy}, first we note that for the initial data
term
\begin{equation}
  \nonumber
   \| W_n^\init \|_\H \to \| W^\init \|_\H\,.
\end{equation}
The convergence of $W_n$ in $C([0,\infty);\wH)$ and
$L^2_{loc}([0,\infty);\H)$, together with the fact that the norm of
the weak limit is an eventual lower bound to the norms of the
sequence, yields
\begin{equation}
   \nonumber
   \| W \|_\H^{\,2}
   \leq \liminf_{n\to\infty} \| W_n(t) \|_\H^{\,2} \,.
\end{equation}
Similarly, the convergence of $W_n$ in
$\mbox{w-}L^2_{loc}([0,\infty); \wV)$ implies
\begin{equation}
  \nonumber
  - \int_0^t (W \,|\, \DDBar W)_{\HBbb} \, \dt
  \leq - \liminf_{n\to\infty}
         \int_0^t (W_n \,|\, \DDBar W_n)_{\HBbb} \, \dt \,.
\end{equation}
Thus, we finish the proof of global Leray type weak solutions. \qed

\subsection{Uniqueness Theorem}
Uniqueness can never be asserted by such a compactness argument, but
generally requires the knowledge of additional regularity of the
solution.  For example, here we will prove the following weak-strong
theorem.

\begin{Thm}
  \label{GWS-Thm-4}
Let $U_1,U_2\in C([0,\infty);\wH)\cap L^2_\loc(\dt;\VBbb)$ be two weak
solutions of the WND system with initial data $U_1^\init,U_2^\init\in\H$.
Let $s>\max\{d/2,1\}$.  If $U_1\in L^2([0,T];\H^s)\cap L^1([0,T];\VBbb^s)$
for some $T>0$ then $U_1\in C([0,T];\wH^{s-1})$ and for every $t\in[0,T]$
one has the energy equality
\begin{equation}
  \label{GWS-40}
  \tfrac12 \big\| U_1(t) \big\|_\H^{\,2}
  - \int_0^t \big( U_1 \,|\, \DDBar U_1 \big)_\H \, \dt'
  = \tfrac12 \big\| U_1^\init \big\|_\H^{\,2} \,,
\end{equation}
and the stability bound
\begin{equation}
  \label{GWS-41}
  \big\| U_2(t) - U_1(t) \big\|_\H
  \leq \exp\!\left( C_\QQ^s \int_0^t \big\| \GRAD U_1(t) \big\|_{\H^s}
                    \, \dt' \right)
       \big\| U_2^\init - U_1^\init \big\|_\H \,.
\end{equation}
In particular, if $U_2^\init=U_1^\init$ then $U_2(t)=U_1(t)$ for every
$t\in[0,T]$.
\end{Thm}

\smallskip
\noindent {\bf Remark.}  Equation (\ref{GWS-40}) is simply the
assertion that the energy inequality satisfied by the strong
solution $U_1$ is in fact an equality.  The bound (\ref{GWS-41}) is
a basic weak-strong stability bound, from which the uniqueness
assertion follows immediately.

\smallskip

   The key to the proof of Theorem \ref{GWS-Thm-4} will be provided by
the following lemma.

\begin{Lem}
  \label{GWS-Lem-4}
Let $U_1,U_2\in C([0,\infty);\wH)\cap L^2_\loc(\dt;\VBbb)$ be two weak
solutions of the WND system with initial data $U_1^\init,U_2^\init\in\H$.
Let $s>\max\{d/2,1\}$.  If $U_1\in L^2([0,T];\H^s)\cap L^1([0,T];\VBbb^s)$
for some $T>0$ then $U_1\in C([0,T];\wH^{s-1})$ and for every $t\in[0,T]$
one has
\begin{equation}
  \label{GWS-42}
\begin{aligned}
  \big( U_1(t) \,|\, U_2(t) \big)_\H
  + \int_0^t \big( \QQBar(U_1,U_1) \,|\, U_2 \big)_\H
             + \big( U_1 \,|\, \QQBar(U_2,U_2) \big)_\H \, \dt' &
\\
  - \int_0^t \big( \DDBar U_1 \,|\, U_2 \big)_\H
             + \big( U_1 \,|\, \DDBar U_2 \big)_\H \, \dt'
  & = \big( U_1^\init \,|\, U_2^\init \big)_\H \,.
\end{aligned}
\end{equation}
\end{Lem}

\smallskip
\noindent
This lemma will be proved later.  Now we will use it to prove
Theorem \ref{GWS-Thm-4}.

\smallskip
\noindent
{\bf Proof of Theorem \ref{GWS-Thm-4}.}  The energy equality
(\ref{GWS-40}) follows by setting $U_2=U_1$ in equation (\ref{GWS-42})
of Lemma \ref{GWS-Lem-4}, using the cyclic identity (\ref{AOp-9})
to see that $\big(U_1\,|\,\QQBar(U_1,U_1)\big)_\H=0$, and multiplying
the result by $\frac12$.

   We now derive the stability bound (\ref{GWS-41}).  Add the energy
inequalities for $U_1$ and $U_2$ and subtract equation (\ref{GWS-42})
from the result to obtain
\begin{equation}
  \label{GWS-44}
\begin{aligned}
  \tfrac12 \big\| U_2(t) - U_1(t) \big\|_\H^{\,2}
  - \int_0^t \big( \QQBar(U_1,U_1) \,|\, U_2 \big)_\H
             + \big( U_1 \,|\, \QQBar(U_2,U_2) \big)_\H \, \dt' &
\\
  - \int_0^t \big( (U_2 - U_1) \,|\,
                   \DDBar (U_2 - U_1) \big)_\H \, \dt'
  & \leq \tfrac12 \big\| U_2^\init - U_1^\init \big\|_\H^{\,2} \,.
\end{aligned}
\end{equation}
Upon letting $W=U_2-U_1$ (so that $U_2=U_1+W$), we see that
\begin{equation}
  \label{GWS-45}
\begin{aligned}
  \big( \QQBar(U_1,U_1) \,|\, U_2 \big)_\H
  + \big( U_1 \,|\, \QQBar(U_2,U_2) \big)_\H
  & = \big( U_1 \,|\, \QQBar(U_1,U_1) \big)_\H
      + \big( W \,|\, \QQBar(U_1,U_1) \big)_\H
\\
  & \quad \,
      + \big( U_1 \,|\, \QQBar(U_1,U_1) \big)_\H
      + \big( U_1 \,|\, \QQBar(W,U_1) \big)_\H
\\
  & \quad \,
      + \big( U_1 \,|\, \QQBar(U_1,W) \big)_\H
      + \big( U_1 \,|\, \QQBar(W,W) \big)_\H \,.
\end{aligned}
\end{equation}
The cyclic identity (\ref{AOp-9}) implies that
$\big(U_1\,|\,\QQBar(U_1,U_1)\big)_\H=0$ and
\begin{equation}
  \nonumber
  \big( W \,|\, \QQBar(U_1,U_1) \big)_\H
  + \big( U_1 \,|\, \QQBar(W,U_1) \big)_\H
  + \big( U_1 \,|\, \QQBar(U_1,W) \big)_\H = 0 \,.
\end{equation}
We thereby see that relation (\ref{GWS-45}) reduces to
\begin{equation}
  \nonumber
  \big( \QQBar(U_1,U_1) \,|\, U_2 \big)_\H
  + \big( U_1 \,|\, \QQBar(U_2,U_2) \big)_\H
  = \big( U_1 \,|\, \QQBar(W,W) \big)_\H \,.
\end{equation}
When this relation is placed into (\ref{GWS-44}), we obtain
\begin{equation}
  \label{GWS-46}
  \tfrac12 \big\| W(t) \big\|_\H^{\,2}
  - \int_0^t \big( U_1 \,|\, \QQBar(W,W) \big)_\H \, \dt'
  - \int_0^t \big( W \,|\, \DDBar W \big)_\H \, \dt'
  \leq \tfrac12 \big\| W^\init \big\|_\H^{\,2} \,,
\end{equation}
where $W^\init=U_2^\init-U_1^\init$.  The third bound in
(\ref{AOp-4}) gives
\begin{equation}
  \nonumber
  \big| \big( U_1 \,|\, \QQBar(W,W) \big)_\H \big|
  \leq C_\QQ^s \big\| \GRAD U_1 \big\|_{\H^s}
               \big\| W \big\|_\H^{\,2} \,.
\end{equation}
We combine this bound with the fact
$-\big( W\,|\,\DDBar W\big)_\H\geq0$ to see that (\ref{GWS-46})
yields the inequality
\begin{equation}
  \nonumber
  \tfrac12 \big\| W(t) \big\|_\H^{\,2}
  \leq \tfrac12 \big\| W^\init \big\|_\H^{\,2}
       + \int_0^t C_\QQ^s \big\| \GRAD U_1 \big\|_{\H^s}
                  \big\| W \big\|_\H^{\,2} \, \dt' \,.
\end{equation}
The stability bound (\ref{GWS-41}) then follows by the Gronwall
Lemma. \qed

\smallskip

   All that remains is to prove Lemma \ref{GWS-Lem-4}.

\smallskip
\noindent
{\bf Proof of Lemma \ref{GWS-Lem-4}.}  We begin by showing that
$U_1\in L^2([0,T];\H^s)\cap L^1([0,T];\V^s)$ implies
$U_1\in C([0,T];\wH^{s-1})$.  From equation \eqref{GWS-1}
we see that for every $W \in \H^{-(s-1)}$ we have
\begin{equation}
  \nonumber
  \big( W \,|\, U_1(t_2) - U_1(t_1) \big)_\H
  = - \int^{t_2}_{t_1}
          \big( W \,|\, \AA U_1(t) + \QQBar(U_1(t), U_1(t))
                        - \DDBar U_1(t) \big)_\H \, \dt \,.
\end{equation}
We have the following estimates:
\begin{equation}
  \nonumber
  \left| \int^{t_2}_{t_1}
             \big( W \,|\, \AA U_1(t) \big)_\H \, \dt \right|
  \leq C_\AA \, \| W \|_{\H^{-(s-1)}}
       \int_{t_1}^{t_2} \| U_1(t) \|_{\H^s} \, \dt \,.
\end{equation}
By the first bound in (\ref{AOp-3}) of Lemma \ref{AOp-Lem-2}
with $s'=s$ we have
\begin{equation}
  \nonumber
\begin{aligned}
  \left| \int^{t_2}_{t_1}
             \big( W \,|\, \DDBar U_1(t) \big)_\H \, \dt \right|
  & \leq C_\DD \, \| \GRAD W \|_{\H^{-s}}
         \int^{t_2}_{t_1} \| \GRAD U_1(t) \|_{\H^s} \, \dt
\\
  & \leq C_\DD \, \| W \|_{\H^{-(s-1)}}
         \int^{t_2}_{t_1} \| U_1(t) \|_{\V^s} \, \dt \,.
\end{aligned}
\end{equation}
By the third bound in (\ref{AOp-3}) of Lemma \ref{AOp-Lem-2}
with $s'=s-1$ we have
\begin{equation}
  \nonumber
  \left| \int^{t_2}_{t_1}
             \big( W \,|\, \QQBar(U_1(t),U_1(t)) \big)_\H \, \dt
         \right|
  \leq 2^s C_\QQ^{s-1} \, \| W \|_{\H^{-(s-1)}}
       \int^{t_2}_{t_1} \| U_1(t) \|_{\H^s}^{\,2} \, \dt \,.
\end{equation}
Because $U_1\in L^2([0,T];\H^s)\cap L^1([0,T];\V^s)$ it follows
that $(W\,|\,U_1(t))_\H$ is a continuous function of $t$ over
$[0,T]$ for every $W\in\H^{s-1}$.  Hence, $U_1 \in C([0,T];\wH^{s-1})$.

   We now prove that (\ref{GWS-42}) holds for every $t\in[0,T]$.  Let
$\alpha\in\DD(\Rd)$ and $\beta\in\DD(\R)$ be mollifiers such that
$\alpha\geq0$, $\beta\geq0$, supp$(\beta)\subset(-\infty,0]$ and
\begin{equation}
  \nonumber
  \int_\Rd \alpha(x) \, \dx = 1 \,, \qquad
  \int_\R \beta(t) \, \dt = 1 \,.
\end{equation}
For every $\eps>0$ define $\Theta_\eps\in \DD(\Td\times\R)$ by
\begin{equation}
  \nonumber
  \Theta_\eps(x,t)
  = \frac{1}{\eps^{d+1}} \,
    \sum_{l\in\Zd} \alpha\!\bigg( \frac{x + 2 \pi l}{\eps} \bigg) \,
    \beta\!\bigg( \frac{t}{\eps} \bigg) \,.
\end{equation}
For each $i=1,2$ and $\eps>0$ define $U_{i\eps}=\Theta_\eps*U_i$, so
that for every $x\in\Td$ and $t\geq0$ we have
\begin{equation}
  \nonumber
  U_{i\eps}(x,t)
  = \big( \Theta_\eps * U_i \big)(x,t)
  = \int_0^\infty \!\! \int_\Td
        \Theta_\eps(x - x', t - t') \,
        U_i(x',t') \, \dx' \, \dt' \,.
\end{equation}

   Each $U_{i\eps}$ is a smooth function over $\Rd\times\R_+$ that
satisfies
\begin{equation}
  \nonumber
  \del_t U_{i\eps} + \AA U_{i\eps} + \Theta_\eps * \QQBar(U_i,U_i)
  = \DDBar  U_{i\eps} \,.
\end{equation}
Because $\AA$ is skew-adjoint we thereby see that
\begin{equation}
  \nonumber
\begin{aligned}
  \frac{\dee}{\dt} \big( U_{1\eps} \,|\, U_{2\eps} \big)_\H
  & = \big( \del_t U_{1\eps} \,|\, U_{2\eps} \big)_\H
      + \big( U_{1\eps} \,|\, \del_t U_{2\eps} \big)_\H
\\
  & = - \big( \Theta_\eps * \QQBar(U_1,U_1) \,|\, U_{2\eps} \big)_\H
      - \big( U_{1\eps} \,|\, \Theta_\eps * \QQBar(U_2,U_2) \big)_\H
\\
  & \quad \,
      + \big( \DDBar U_{1\eps} \,|\, U_{2\eps} \big)_\H
      + \big( U_{1\eps} \,|\, \DDBar U_{2\eps} \big)_\H \,.
\end{aligned}
\end{equation}
Upon integrating this equation over $[0,t]$ we obtain
\begin{equation}
  \label{GWS-47}
\begin{aligned}
  \big( U_{1\eps}(t) \,|\, U_{2\eps}(t) \big)_\H
  & + \int_0^t \big( \Theta_\eps * \QQBar(U_1,U_1) \,|\,
                     U_{2\eps} \big)_\H
             + \big( U_{1\eps} \,|\,
                     \Theta_\eps * \QQBar(U_2,U_2) \big)_\H \, \dt'
\\
  & - \int_0^t \big( \DDBar U_{1\eps} \,|\, U_{2\eps} \big)_\H
             + \big( U_{1\eps} \,|\, \DDBar U_{2\eps} \big)_\H \, \dt'
\\
  & = \big( U_{1\eps}(0) \,|\, U_{2\eps}(0) \big)_\H \,.
\end{aligned}
\end{equation}

   We claim that
\begin{equation}
  \label{GWS-48}
\begin{aligned}
  \lim_{\eps \to 0} \big( U_{1\eps}(t) \,|\, U_{2\eps}(t) \big)_\H
  & = \big( U_1(t) \,|\, U_2(t) \big)_\H \,,
\\
  \lim_{\eps \to 0} \big( U_{1\eps}(0) \,|\, U_{2\eps}(0) \big)_\H
  & = \big( U_1^\init \,|\, U_2^\init \big)_\H \,,
\\
  \lim_{\eps \to 0}
      \int_0^t \big( \Theta_\eps * \QQBar(U_1,U_1) \,|\,
                     U_{2\eps} \big)_\H \, \dt'
  & = \int_0^t \big( \QQBar(U_1,U_1) \,|\, U_2 \big)_\H \, \dt' \,,
\\
  \lim_{\eps \to 0}
      \int_0^t \big( U_{1\eps} \,|\,
                     \Theta_\eps * \QQBar(U_2,U_2) \big)_\H \, \dt'
  & = \int_0^t \big( U_1 \,|\, \QQBar(U_2,U_2) \big)_\H \, \dt' \,,
\\
  \lim_{\eps \to 0}
      \int_0^t \big( \DDBar U_{1\eps} \,|\, U_{2\eps} \big)_\H
             + \big( U_{1\eps} \,|\, \DDBar U_{2\eps} \big)_\H \, \dt'
  & = \int_0^t \big( \DDBar U_1 \,|\, U_2 \big)_\H
             + \big( U_1 \,|\, \DDBar U_2 \big)_\H \, \dt' \,.
\end{aligned}
\end{equation}
Once these limits are established we can then pass to the limit in
(\ref{GWS-47}) to obtain (\ref{GWS-42}) and thereby complete the
proof of Lemma \ref{GWS-Lem-4}.  The limits (\ref{GWS-48}) are
established by using the bounds from Lemma \ref{AOp-Lem-2} and
the convergence and boundedness properties of convolution.

   Now the {\em first two} limits are direct consequence of the
convergence property of convolution in $C([0,T];\wH^{s-1})$.
Because $U_1$ and $U_2$ are continuous in time in $\wH^{s-1}$
and $\H$ respectively, we have, for any $t > 0$,
\begin{equation}
  \nonumber
\begin{aligned}
  & U_{1\eps}(t) \to U_1(t) \quad \mbox{in}\quad\!
  \wH^{(s-1)}\quad\mbox{thus}\quad U_{1\eps}(t) \to
U_1(t)\quad \mbox{in}\quad\! \H\,, \\
 & U_{2\eps}(t) \to U_2(t)\quad \mbox{in}\quad \! \wH\,,
\end{aligned}
\end{equation}
as $\eps \to 0$, which imply that for every $t > 0$,
\begin{equation}
  \nonumber
  \big( U_{1\eps}(t) \,|\, U_{2\eps}(t) \big)_\H
  \to \big( U_1(t) \,|\, U_2(t) \big)_\H \,,
\end{equation}
as $\eps \to 0$. Thus we prove the first two limits.

  To prove the {\em third} limit, first, we have
\begin{equation}
  \label{3rd}
\begin{aligned}
  & \int_0^t \big( \Theta_\eps * \QQBar(U_1,U_1) \,|\,
                   U_{2\eps} \big)_\H
             - \big( \QQBar(U_1,U_1) \,|\, U_2 \big)_\H \, \dt'
\\
  & = \int_0^t \big( \Theta_\eps * \QQBar(U_1,U_1) - \QQBar(U_1,U_1)
                     \,|\, U_{2\eps} \big)_\H \, \dt'
      + \int_0^t \big( \QQBar(U_1,U_1)
                       \,|\, U_{2\eps} - U_2 \big)_\H \, \dt' \,.
 \end{aligned}
\end{equation}
Note that $U_1$ is a strong solution, i.e.
$U_1\in L^\infty([0,T];\H^s)\cap L^2([0,T];\V^s)$ for $s > d/2$,
and the structure of $\QQBar(U_1,U_1)$ is the derivative of $U_1$
multiplying $U_1$, thus
\begin{equation}
  \nonumber
  \QQBar(U_1,U_1) \in L^2([0,T] \,; \H) \,.
\end{equation}
Then the convergence in $L^2([0,T];\H)$ and the boundedness of
mollifier imply that each integral on the right-hand side of
\eqref{3rd} goes to $0$ as $\eps\to 0$.

   We leave the proof of the fourth limit to the last step because
it is the hardest one.  We prove the {\em fifth} limit first.
Applying the first inequality in Lemma \ref{AOp-Lem-2} with $s'=0$
and the H\"older inequality, we have
\begin{equation}
  \nonumber
\begin{aligned}
  \left| \int_0^t \big( \DDBar U_{1\eps} \,|\, U_{2\eps} \big)_\H
                  - \big( \DDBar U_1 \,|\, U_2 \big)_\H
                  \, \dt' \right|
  & = \left| \int_0^t \big( \DDBar U_{1\eps} - \DDBar U_1
                            \,|\, U_{2\eps} \big)_\H
                      + \big( \DDBar U_1 \,|\, U_{2\eps} - U_2
                              \big)_\H \, \dt' \right|
\\
  & \leq C_\DD \| \GRAD (U_{1\eps} - U_1) \|_{L^2([0,T] \,; \H)}
               \| \GRAD U_{2\eps} \|_{L^2([0,T] \,; \H)}
\\
  & \quad \,
         + C_\DD \| \GRAD U_1 \|_{L^2([0,T] \,; \H)}
           \| \GRAD (U_{2\eps} - U_2) \|_{L^2([0,T] \,; \H)} \,.
\end{aligned}
\end{equation}
Note that $\GRAD U_{i\eps}=\big( \GRAD U_i \big)_\eps\,,\quad
\mbox{for}\quad\! i=1,2\,,$ the boundedness of $\GRAD U_{i\eps}$ in
$L^2\big( [0,T]\,; \H \big)$,  and the convergence
\begin{equation}
  \nonumber
  \big( \GRAD U_i \big)_\eps \to \GRAD U_i \quad
  \mbox{in} \quad\! L^2\big( [0,T]\,; \H \big) \quad
  \mbox{as} \quad\! \eps \to 0\,,
\end{equation}
we finish the proof of the fifth limit.

   The main difficulty is that $U_2$ is only a Leray weak solution,
so $\QQBar(U_2,U_2)$ is not in $L^2$, thus the method to prove the
third limit is not applicable here. However, we have the following
identity: for any functions $U, V$ so that $U^T G_o V  \in
L^1([0,T]\,; L^1(\dx))$ and $T>0$,
\begin{equation}
  \label{changing-order}
  \int_0^T\big( \Theta_\eps * U\,|\, V \big)_\H \,\dt
  = \int_0^T\big( U \,|\, \Theta_\eps * V \big)_\H \,\dt \,.
\end{equation}
{\bf Proof of \eqref{changing-order}}: By changing the order of
integration,
\begin{equation}
  \nonumber
\begin{aligned}
  & \int_0^T \big( \Theta_\eps * U \,|\, V \big)_\H \,\dt
\\
  & = \int_0^T\int_{\T^d} \int_0^\infty\int_{\T^d}
         \tfrac{1}{\eps^{d+1}}
         \alpha\big( \tfrac{x-x'}{\eps} \big)
         \beta\big( \tfrac{t-t'}{\eps} \big) \, U^T(x',t') \, \dx' \, \dt'
         \, G_o\,V(x,t)\,\dx\, \dt\\
  & = \int_0^T \int_{\T^d} \int^T_{t'} \int_{\T^d}
          \tfrac{1}{\eps^{d+1}}
          \alpha\big( \tfrac{x-x'}{\eps} \big)
          \beta\big( \tfrac{t-t'}{\eps} \big) \, V^T(x,t) \, \dx \, \dt
          \, G_o \,U(x',t') \, \dx' \, \dt'
\\
  & = \int_0^T \big( U \,|\, \Theta_\eps * V \big)_\H \,\dt \,.
\end{aligned}
\end{equation}

   Applying \eqref{changing-order} to $U_1$ and $\QQBar(U_2,U_2)$, we
have
\begin{equation}
  \nonumber
  \int_0^t \big( U_{1\eps} \,|\,
                 \Theta_\eps * \QQBar(U_2,U_2) \big)_\H \, \dt'
  = \int_0^t \big( \Theta_\eps * U_{1\eps} \,|\,
                   \QQBar(U_2,U_2) \big)_\H \, \dt' \,.
\end{equation}
Thus from the second inequality in Lemma \ref{AOp-Lem-2},
\begin{equation}
  \nonumber
\begin{aligned}
  & \left| \int_0^t \big( \Theta_\eps * U_{1\eps}-U_1 \,|\,
                          \QQBar(U_2,U_2) \big)_\H \, \dt' \right|
\\
  & \leq C^s_\QQ T^\frac12
         \int_0^t \| \GRAD (\Theta_\eps
                   * U_{1\eps} - U_1) \|_{L^\infty([0,T];\H^s)} \, \dt'
         \| U_2 \|^2_{L^\infty(\dt;\H)} \,.
\end{aligned}
\end{equation}
Note that
\begin{equation}\nonumber
  \GRAD \big( \Theta_\eps * U_{\eps} \big)
  = \Theta_\eps * (\GRAD U)_\eps \,.
\end{equation}
We claim that for any $W \in L^2([0, T];\H^s)$,
\begin{equation}
  \nonumber
  \lim_{\eps \to 0}
      \| \Theta_\eps * W_\eps - W \|_{L^2([0,T]\,; \H^s)}
  \to 0 \quad \mbox{as $\eps \to 0$} \,,
\end{equation}
which simply followed from the triangle inequality
\begin{equation}
  \nonumber
\begin{aligned}
  \| \Theta_\eps * W_\eps - W \|_{L^2([0,T];\H^s)}
  \leq \|\Theta_\eps * W_\eps - W_\eps\|_{L^2([0,T];\H^s)}
       + \| W_\eps - W \|_{L^2([0,T];\H^s)} \,,
\end{aligned}
\end{equation}
and the convergence property of convolution in $L^2([0,T];\H^s)$.  We
thereby prove the limit.

   We have now established all the limits asserted in \eqref{GWS-48}
and have thereby completed the proof of Lemma \ref{GWS-Lem-4}. \qed

%%%%%%%%%%%%%%%%%%%%%%%%%%%%%%%%%%%%%%%%%%%%%%%%%%%%%%%%%%%%%%%%%%%%%%
%%  Section 7: Application to the Compressible Navier-Stokes System %%
%%%%%%%%%%%%%%%%%%%%%%%%%%%%%%%%%%%%%%%%%%%%%%%%%%%%%%%%%%%%%%%%%%%%%%

\section{Application to the Compressible Navier-Stokes System}
\label{CNS}

\subsection{Compressible Navier-Stokes System}

   The compressible Navier-Stokes system of gas dynamics is an
important example of a nonsingular hyperbolic-parabolic system with
the strictly convex entropy.  It governs the mass density $\rho(x,t)$,
bulk velocity $u(x,t)$, and temperature $\theta(x,t)$ over
$\Omega\subset\Rd$ in the form
\begin{equation}
  \label{CNS-1}
\begin{aligned}
  \del_t \rho + \DIV (\rho u) & = 0 \,,
\\
  \del_t (\rho u)
  + \DIV \big( \rho u \otimes u + p I + S \big) & = 0 \,,
\\
  \del_t (\tfrac12 \rho |u|^2 + \rho \vareps)
  + \DIV \big (\tfrac12 \rho |u|^2 u + \rho \vareps u + p u
                                     + S \DOT u + q \big) & = 0 \,,
\end{aligned}
\end{equation}
where the specific energy $\vareps$ and pressure $p$ are given by
thermodynamic equations-of-state, $\vareps=\vareps(\rho,\theta)$ and
$p=p(\rho,\theta)$, while the stress $S$ and heat flux $q$ are given
by the constitutive relations
\begin{equation}
  \label{CNS-2}
  S = - \mu \, \big( \GRAD u + (\GRAD u)^T
                     - \tfrac2{\D} I \, \DIV u \big)
        - \lambda \, I \, \DIV u \,, \qquad
  q = - \kappa \, \GRAD \theta \,.
\end{equation}
Here the coefficients of shear viscosity $\mu$, bulk viscosity
$\lambda$ and thermal conductivity $\kappa$ are given by formulas
$\mu=\mu(\rho,\theta)$, $\lambda=\lambda(\rho,\theta)$, and
$\kappa=\kappa(\rho,\theta)$ that come either from a nonequilibrium
(kinetic) theory or from fits to experimental data, while $\D$ is
the dimension of the underlying microscopic world --- usually $\D=3$.
We require that $\D\geq\max\{2,d\}$.

   Equations (\ref{CNS-1}) express the local conservation of mass,
momentum, and energy.  The constitutive relations (\ref{CNS-2}) for
$S$ and $q$ model viscosity and thermal conductivity, which arise due
to deviations of the gas from local thermodynamic equilibrium.
Equations (\ref{CNS-1}) reduce to the compressible Euler system when
one sets $S=0$ and $q=0$.

   The thermodynamic equations-of-state for the specific energy and
pressure, $\vareps=\vareps(\rho,\theta)$ and $p=p(\rho,\theta)$, are
assumed to be twice continuously differentiable over
$(\rho,\theta)\in\R_+^2$ and to satisfy
\begin{equation}
  \label{CNS-3}
  \del_\theta \vareps(\rho,\theta) > 0 \,, \qquad
  \del_\rho p(\rho,\theta) > 0 \,, \qquad
  \hbox{for every $\rho>0$ and $\theta>0$} \,.
\end{equation}
In addition, they must satisfy the Maxwell relation
\begin{equation*}
  \rho^2 \del_\rho \vareps
  + \theta^2 \del_\theta \bigg( \frac{p}{\theta} \bigg) = 0 \,.
\end{equation*}
This implies the existence of a function $\sig=\sig(\rho,\theta)$
that satisfies the differential relation
\begin{equation*}
  \dee \bigg( \sig - \frac{\vareps}{\theta} \bigg)
  = - \vareps \, \dee \bigg( \frac{1}{\theta} \bigg)
    + \frac{p}{\theta} \, \dee \bigg( \frac{1}{\rho} \bigg) \,.
\end{equation*}
This is equivalent to
\begin{equation}
  \label{CNS-4}
  \dee \sig
  = \frac{1}{\theta} \, \dee \vareps
    - \frac{p}{\rho^2 \theta} \, \dee \rho \,.
\end{equation}
We can identify $\sig$ with the specific entropy.

   We make use of the convective form of system (\ref{CNS-1}),
\begin{equation}
  \label{CNS-5}
\begin{aligned}
  \del_t \rho + u \DOT \GRAD \rho + \rho \DIV u & = 0 \,, \Big.
\\
  \rho (\del_t u + u \DOT \GRAD u)
  + \GRAD p + \DIV S & = 0 \,, \Big.
\\
  \rho (\del_t \vareps + u \DOT \GRAD \vareps)
  + p \DIV u + S \DDOT \GRAD u + \DIV q & = 0 \,, \Big.
\end{aligned}
\end{equation}
and the differential specific entropy relation (\ref{CNS-4}) to see
that
\begin{equation}
  \nonumber
\begin{aligned}
  \rho (\del_t \sig + u \DOT \GRAD \sig)
  & = \frac{\rho}{\theta} \,
      (\del_t \vareps + u \DOT \GRAD \vareps)
      - \frac{p}{\rho \theta} \,
      (\del_t \rho + u \DOT \GRAD \rho)
    = - \frac{1}{\theta} \, S \DDOT \GRAD u
      - \frac{1}{\theta} \,  \DIV q \,.
\end{aligned}
\end{equation}
This can be put into the divergence form
\begin{equation}
  \label{CNS-6}
  \del_t (\rho \sig)
  + \DIV \bigg( \rho u \sig + \frac{q}{\theta} \bigg)
  = - \frac{1}{\theta} \, S \DDOT \GRAD u
    - \frac{1}{\theta^2} \, q \DOT \GRAD \theta \,.
\end{equation}
The local form of the second law of thermodynamics and the
constitutive relations (\ref{CNS-2}) imply that for any values of
$\rho$, $\theta$, $\GRAD u+(\GRAD u)^T-\tfrac2{\D}I\DIV u$, $\DIV u$,
and $\GRAD \theta$ one has the inequality
\begin{equation}
  \nonumber
\begin{aligned}
   - \frac{1}{\theta} \, S \DDOT \GRAD u
   - \frac{1}{\theta^2} \, q \DOT \GRAD \theta
   & = \frac{\mu}{2} \, \big| \GRAD u + (\GRAD u)^T
                              - \tfrac2{\D} I \, \DIV u \big|^2
        + \lambda \, |\DIV u|^2
        + \kappa \, |\GRAD \theta|^2
     \geq 0 \,.
\end{aligned}
\end{equation}
Because these values can be independently specified at any point in
$\Omega$, the above inequality implies that $\mu(\rho,\theta)\geq0$,
$\lambda(\rho,\theta)\geq0$, and $\kappa(\rho,\theta)\geq0$ for every
$\rho>0$ and $\theta>0$.  Of course, these thermodynamic constraints
are satisfied when $\mu=\lambda=\kappa=0$, which is the case of the
compressible Euler system.  For the compressible Navier-Stokes system
we require that
\begin{equation}
  \label{CNS-7}
  \mu(\rho,\theta) > 0 \,, \quad
  \lambda(\rho,\theta) \geq 0 \,, \quad
  \kappa(\rho,\theta) > 0 \,, \quad
  \hbox{for every $\rho>0$ and $\theta>0$} \,,
\end{equation}
which are also consistent with the thermodynamic constraints.

   Now consider the system over a periodic box $\Omega=\Td$.
Integrating the divergence form of the entropy equation (\ref{CNS-6})
over $\Td$ yields
\begin{equation}
  \label{CNS-8}
\begin{aligned}
  \frac{\dee}{\dt} \int_\Td \rho \sig \, \dx
  & = \int_\Td
          \bigg[ \frac{\mu}{2} \,
                 \big| \GRAD u + (\GRAD u)^T
                       - \tfrac2{\D} I \, \DIV u \big|^2
                 + \lambda \, |\DIV u|^2
                 + \kappa \, |\GRAD \theta|^2 \bigg] \, \dx \,.
\end{aligned}
\end{equation}
One sees from (\ref{CNS-7}) and (\ref{CNS-8}) that any stationary
classical solution of the system for which $\rho>0$ and $\theta>0$
must satisfy
\begin{equation}
  \nonumber
  \GRAD u + (\GRAD u)^T - \tfrac2{\D} I \, \DIV u = 0 \,, \qquad
  \GRAD \theta = 0 \,.
\end{equation}
Because $u$ is periodic, one can use the first equation above to argue
that $\GRAD u=0$.  It then follows from (\ref{CNS-2}) and
(\ref{CNS-5}) that $\GRAD p(\rho,\theta)=0$, which by the second
equation above and (\ref{CNS-3}) yields $\GRAD\rho=0$.  This shows
that the compressible Navier-Stokes system (\ref{CNS-1}-\ref{CNS-2})
satisfies the nonsingularity condition (\ref{HP-10}).

   The compressible Navier-Stokes system (\ref{CNS-1}-\ref{CNS-2}) is
a nonsingular hyperbolic-parabolic system with a strictly convex
entropy given by $H(U)=-\rho\sig(\rho,\theta)$ where $U$ is related
to $\rho$, $u$, and $\theta$ by
\begin{equation}
  \nonumber
  U = \begin{pmatrix}
         \rho  &  \rho u_1 & \cdots & \rho u_\D &
         \tfrac12 \rho |u|^2 + \rho \vareps(\rho,\theta)
      \end{pmatrix}^T \,.
\end{equation}
The set $\UU$ is the range of this mapping restricted to the domain
$\rho>0$, $u\in\Rd$, and $\theta>0$.  Whenever
$\theta\mapsto\vareps(\rho,\theta)$ is a strictly increasing function
from $\R_+$ {\em onto} $\,\R_+$ then $\UU$ is given by
\begin{equation}
  \nonumber
  \UU = \Big\{ \begin{pmatrix}
                  U_0 & U_1 & \cdots & U_\D & U_{\D+1}
               \end{pmatrix}^T \in \R^{\D+2} \,:\, U_0 > 0 \,, \
               2 U_0 U_{\D+1} > U_1^{\,2} + \cdots + U_\D^{\,2}
               \Big\} \,.
\end{equation}
The function $H:\UU\to\R$ will be strictly convex if and only if
$\del_\theta \vareps(\rho,\theta)>0$ and the sound speed is defined
\cite{HLLM}.  These conditions are satisfied by all thermodynamic
equations-of-state that satisfy (\ref{CNS-3}).

   Finally, it is easily checked that the Navier-Stokes system
(\ref{CNS-1}-\ref{CNS-2}) satisfies our Kawashima-type criterion
(\ref{SDA-1}), whereby it also satisfies the Kawashima condition
(\ref{WND-7}).  Moreover, because its weakly nonlinear-dissipative
approximation will be strictly dissipative, our theory of global
weak solutions applies to it.

\subsection{Weakly Compressible Navier-Stokes System}
\label{WCNS}

   A fluid dynamical system that formally includes both the acoustic
and the Stokes systems is the so-called
{\em weakly compressible Stokes system}
\begin{equation*}
\begin{aligned}
  \del_t \rhoTld + \rho_o \DIV \uTld & = 0 \,,
\\
  \rho_o \del_t \uTld
  + \GRAD \big( (\del_\rho p)_o \, \rhoTld
                + (\del_\theta p)_o \, \thetaTld \big)
  & = \mu_o \DIV \big[ \GRAD \uTld + (\GRAD \uTld)^T
                       - \tfrac{2}{\D} \DIV \uTld I \big]
      + \lambda_o \GRAD (\DIV \uTld) \,,
\\
  \rho_o C_o^V \del_t \thetaTld
  + \theta_o (\del_\theta p)_o \, \DIV \uTld
  & = \kappa_o \LAP \thetaTld \,,
\end{aligned}
\end{equation*}
where $C_o^V=(\del_\theta \vareps)_o$ is the specific heat capacity at
constant volume.  A so-called {\em weakly compressible Navier-Stokes
system} that formally includes both the acoustic and the Navier-Stokes
systems is the weakly nonlinear-dissipative approximation.  It
decomposes $\UTld$ into a component governed by the incompressible
Navier-Stokes system and a component governed by a nonlocal quadratic
acoustic equation that couples to the incompressible component.

   The incompressible component $\UTld_{in}=(\rhoTld,\uTld,\thetaTld)$
is governed by
\begin{equation}
  \label{WCNS-1}
\begin{aligned}
  \rho_o \big( \del_t \uTld + \uTld \DOT \GRAD \uTld \big) + \GRAD \pTld
  & = \mu_o \LAP \uTld \,,
\\
  \rho_o C_o^P \big( \del_t \thetaTld + \uTld \DOT \GRAD \thetaTld \big)
  & = \kappa_o \LAP \thetaTld \,,
\end{aligned}
\end{equation}
where
\begin{equation}
  \label{WCNS-2}
  \DIV \uTld = 0 \,, \qquad
  (\del_\rho p)_o \, \rhoTld + (\del_\theta p)_o \, \thetaTld = 0 \,,
\end{equation}
and $C_o^P$ denotes the specific heat capacity at constant pressure,
while the acoustic component $\UTld_{ac}=(\etaTld,\vTld,\chiTld)$ is
governed by
\begin{equation}
  \label{WCNS-3}
  \del_t \UTld_{ac} + \AA \UTld_{ac}
  + \DIV \QBar_{in}(\UTld_{in},\UTld_{ac})
  + \DIV \QBar_{ac}(\UTld_{ac},\UTld_{ac})
  = \nuBar \LAP \UTld_{ac} \,,
\end{equation}
where
\begin{equation}
  \label{WCNS-4}
  \ROT \, \vTld = 0 \,, \qquad
  \chiTld = \frac{\theta_o (\del_\theta p)_o}{\rho_o^2 C_o^V} \,
            \etaTld \,, \qquad \,.
\end{equation}
The diffusion coefficient $\nuBar$ in \eqref{WCNS-4} is given by
\begin{equation}
  \label{Diffusion}
  \nuBar = \frac{2 \frac{\D - 1}{\D} \mu_o + \lambda_o}{2 \rho_o}
           + \frac{\kappa_o}{2 \rho_o C_o^V}
             \frac{\theta_o (\del_\theta p)_o^2}{\rho_o^2 C_o^V c_o^2} \,,
\end{equation}
where $c_o$ is the sound speed defined by
\begin{equation}
   \label{sound-speed}
   c_o^2 = (\del_\rho p)_o
           + \frac{\theta_o (\del_\theta p)_o^2}{\rho_o^2 C_o^V} \,.
\end{equation}
Notice that $\nuBar$ is positive if either $\mu_o$, $\lambda_o$, or $\kappa_o$
is positive! $\QBar_{in}$ and $\QBar_{ac}$ are nonlocal operators defined in the
following way.  We first introduce the orthonomal basis of the acoustic mode
Null$(\AA)^\perp$:
\begin{equation}
  H_k^\pm( x)
  = \sqrt{\tfrac{\theta_o}{2\rho_o}}
    \begin{pmatrix}
       \frac{\rho_o}{ c_o} \\
       \pm\frac{ k}{| k|} \\
       \frac{\theta_o (\del_\theta p)_o}{\rho_o C_o^V c_o}
    \end{pmatrix} e^{i k \DOT x} \,.
\end{equation}
It is easy to check that ${H}^{\pm}_ k$ is the eigenvector of the acoustic
operator $\AA$ with eigenvalues $\pm i c_o|k|$.  Then any acoustic component
$\UTld_{ac}$ can be represented as
\begin{equation}
  \UTld_{ac}= \sum\limits_{k} U^\pm_k {H}^{\pm}_ k( x) \,,
\end{equation}
where $U^\pm_k$ is the coefficient of $\UTld_{ac}$ with respect to the basis
${H}^{\pm}_ k$ under the inner product $(\cdot | \cdot )_{\mathbb{H}}$.  The
nonlocal operator $\QBar_{in}$ can be written as
\begin{equation}
  \QBar_{in}(\UTld_{in},\UTld_{ac})
  = \sum\limits_{\delta, m} \lambda_m^\pm(\UTld_{in}) H^\pm_m(x) \,,
\end{equation}
where
\begin{equation}
  \lambda_m^\pm(\UTld_{in})
  = \sum\limits_{\substack{\pm k \pm l= m \\ |k|=|m|}}
       U^\pm_ k \left[ c_1 \frac{(\widehat{ \tilde{u}}_l\cdot m) k
                                 + (k \cdot m) \widehat{ \tilde{u}}_ l}
                                {| k\| m|}
                       + \frac{\widehat{\tilde{\theta}}_l} {|m|}
                         (c_2 k + c_3 m) \right] \,.
\end{equation}
The nonlocal operator $\QBar_{ac}$ can be written as
\begin{equation}
  \QBar_{in}(\UTld_{ac},\UTld_{ac})
  = \sum\limits_{m} \chi_m^\pm H^\pm_m(x) \,,
\end{equation}
where
\begin{equation}
   \chi^\pm_ m
   = c_4 \sum\limits_{\substack{k+l=m \\
                                \pm(k)|k|+\pm(l)|l|=\pm|m|}}
             U^\pm_ k U^\pm_ l\frac{m}{|m|} \,,
\end{equation}
where the constants $c_1$, $c_2$, $c_3$, and $c_4$ are calculated in
\cite{JL2}.

   We have shown that the weakly compressible Navier-Stokes system has
global weak solutions in $L^2$.  This result includes the Leray
theory, so it cannot be improved easily.  As with the Leray theory,
the key to this result is an ``energy'' dissipation estimate.  We
have, (extending ideas of Masmoudi and Danchin) also adapted a
Littlewood-Payly decomposition to show the acoustic part is unique for
a given incompressible component \cite{JL1}.  Moreover, we have
derived the weakly compressible Navier-Stokes system directly from
the Boltzmann equation for the case $p=\rho\theta$ and
$\vareps=\frac{\D}{2}\theta$ \cite{JL2}.

%%%%%%%%%%%%%%%%%%%%%%%%%%%%%%%%%%%%%%%%%%%%%%%%%%%%%%%%%%%%%%%%%%%%%%


\begin{thebibliography}{99}


\bibitem{BMN97}
A. Babin, A. Mahalov, and B. Nicolaenko,
{\em Global regularity and integrability of 3D Euler and Navier-Stokes
     Equations for Uniformly Rotating Fluids},
Asymptotic Analysis {\bf 15} (1997), 103--150.

\bibitem{BMN99}
A. Babin, A. Mahalov, and B. Nicolaenko,
{\em Global Regularity of 3D Rotating Navier-Stokes Equations
     for Resonant Domains},
Indiana Univ. Math. J. {\bf 48}, (1999) 1133--1176.

\bibitem{BMN01}
A. Babin, A. Mahalov, and B. Nicolaenko,
{\em 3D Navier-Stokes and Euler Equations with Initial Data
     Characterized by Uniformly Large Vorticity},
Indiana Univ. Math. J. {\bf 50} (2001) 1--35.

\bibitem{BHN} S. Bianchini, B. Hanouzet, and R. Natalini,
{\em Asymptotic behavior of smooth solutions for partially dissipative hyperbolic systems with a convex entropy},
 Commun. on Pure \& Appl. Math. {\bf 60} (2007),  no. 11, 1559--1622

\bibitem{CLL}
G.-Q. Chen, C.D. Levermore, and T.-P. Liu, {\em Hyperbolic
Conservation Laws with
     Stiff Relaxation Terms and Entropy},
Commun. Pure \& Appl. Math. {\bf 47} (1994), 787--830.

\bibitem{CF}
P. Constantin and C. Foias,
{\em Navier-Stokes Equations},
Chicago Lectures in Mathematics,
The University of Chicago Press, Chicago, 1988.

\bibitem{Dafermos}
C. Dafermos,
{\em Hyperbolic Conservation Laws in Continuum Physics}, 2nd Edition,
Springer-Verlag, Berlin-Heidelberg, 2005.

\bibitem{Danchin}
R. Danchin,
{\em Zero Mach Number Limit for Compressible Flows
     with Periodic Boundary Conditions},
Amer. J. Math. {\bf 124:6} (2002), 1153--1219.

\bibitem{Desvill}
L. Desvillettes,
{\em Hypocoercivity: The Example of Linear Transport},
%Proceedings of ``UIMP-RSME Santalo Summer School:
%Recent Trends in Partial Differential Equations'', Santander, 2004.
in {\em Recent Trends in Partial Differential Equations},
X. Cabr\'e, J.A. Carrillo, and J.L. V\'azquez eds.,
Contemp. Math. {\bf 409}, 33--54,
Amer. Math. Soc., Providence, RI, 2006

\bibitem{DGL}
C.R. Doering, J.D. Gibbon, and C.D. Levermore,
{\em Weak and Strong solutions of the Complex Ginzburg-Landau Equation},
Physica D {\bf 71} (1994), 285--318.

\bibitem{EM1}
P. Embid and A.J. Majda,
{\em Averaging over Fast Gravity Waves for Geophyical Flows with
     Arbitrary Potenital Vorticity},
Commun. P.D.E. {\bf 21} (1996), 619--658.

\bibitem{EM2}
P. Embid and A.J. Majda,
{\em Low Froude Number Limiting Dynamics for Stably Stratified Flow
     with Small or Finite Rossby Numbers},
Geophys.\& Astrophys. Fluid Dynamics {\bf 87} (1998), 1--50.

\bibitem{FL}
K.O. Friedrichs and P.D. Lax,
{\em Systems of Conservation Equations with a Convex Extension},
Proc. Nat. Acad. Sci. USA {\bf 68} (1971), 1686--1688.

\bibitem{Gallagher1}
I. Gallagher
{\em Asymptotics of the solutions of hyperbolic equations with a
     skew-symmetric perturbation},
Journal of Differential Equations, {\bf 150} (1998), 363--384.

\bibitem{Gallagher2}
I. Gallagher
{\em Applications of Schochet's Methods to Parabolic Equations},
Journal de Mathématiques Pures et Appliquées, {\bf 77} (1998), 989--1054.

\bibitem{Godunov}
S. K. Godunov, {\em An Interesting Class of Quasilinear Systems},
Dokl. Acad. Nauk SSSR {\bf 139} (1961), 521--523.

\bibitem{HLLM}
A. Harten, P.D. Lax, C.D. Levermore, and W.J. Morokoff,
{\em Convex Entropies and Hyperbolicity for General Euler Equations},
SIAM J. Num. Anal. {\bf 35} (1998), 2117--2127.

\bibitem{HZ1}
D. Hoff and K. Zumbrun,
{\em Multi-Dimensional Diffusion Waves for the Navier-Stokes Equations
     of Compressible Flow},
Indiana Univ. Math. J. {\bf 44:2} (1995), 603--676.

\bibitem{HZ2}
D. Hoff and K. Zumbrun,
{\em Pointwise Decay Estimates for Multidimensional Navier-Stokes
     Diffusion Waves},
Z. Angew. Math. Phys. {\bf 48:4} (1997), 597--614.

\bibitem{Hormander}
L. H\"{o}rmander,
{\em Hypoelliptic second order differential operators}, Acta Math. {\bf 119} (1967), 147--171.

\bibitem{JL1}
N. Jiang and C.D. Levermore,
{\em Global Existence and Regularity of
     Weakly Compressible Gas Dynamics},
Preprint 2009.

\bibitem{JL2}
N. Jiang and C.D. Levermore,
{\em Weakly Compressible Stokes Dynamics of the Boltzmann Equation},
Preprint 2009.

\bibitem{Ka}
S. Kawashima,
{\em Systems of Hyperbolic-Parabolic Composite Type,
     with Application to the Equations of Magnetohydrodynamics},
Doctoral Thesis, Kyoto University, 1984.

\bibitem{KaSh}
Y. Shiazuta and S. Kawashima,
{\em System of equations of hyperbolic-parabolic type with applications to the discrete Boltzmann equation},
Tohoku Math. J. {\bf 14} (1985), 249-275.

\bibitem{Kawashima}
S. Kawashima,
{\em Large-time behaviour of solutions to hyperbolic-parabolic systems of conservation laws and applications},  Proc. Roy. Soc. Edinburgh Sect. A  {\bf 106}  (1987),  no. 1-2, 169--194.

\bibitem{Kawashima-BE}
S. Kawashima,
{\em Large-time behavior of solutions of the discrete Boltzmann equation},   Comm. Math. Phys.  {\bf 109}  (1987),  no. 4, 563--589.

\bibitem{La}
O.A. Ladyzhenskaya,
{\em The Mathematical Theory of Viscous Incompressible Flow},
Gordon and Breach, New York, 1969

\bibitem{Le34} J. Leray,
{\em Sur le mouvement d'un fluide visqueux emplissant l'espace},
Acta Math. {\bf 63} (1934), 193--248.

\bibitem{LM1}
P.L. Lions and N. Masmoudi,
{\em Incompressible limit for a viscous compressible fluid},
J. Math. Pures Appli. {\bf 77} (1998), 585--627.

\bibitem{LM2}
P. L. Lions and N. Masmoudi,
{\em Une approche locale de la limite incompressible},
C. R. Acad. Sci. Paris Sr. I Math. {\bf 329} (1999), 387--392.

\bibitem{LZ}
T.-P. Liu and Y. Zeng,
{\em Large Time Behavior of Solutions for General Quasilinear
     Hyperbolic-Parabolic Systems of Conservation Laws},
Mem. Amer. Math. Soc. {\bf 125} (1997), 599.

\bibitem{Masmoudi}
N. Masmoudi,
{\em Incompressible, Inviscid Limit of the
     Compressible Navier-Stokes System},
Ann. Inst. H. Poincar\'{e} {\bf 18} (2001), 199--224.


\bibitem{PuSai04}
M. Puel and L. Saint-Raymond,
{\em Quasineutral Limit for the Relativistic Vlasov-Maxwell System},
Asymptotic Analysis {\bf 40} (2004), 303--352.

\bibitem{Ru-Serre}
T. Ruggeri and D. Serre,
{\em Stability of constant equilibrium state for dissipative balance laws system with a convex entropy},  Quart. Appl. Math.  {\bf 62}  (2004),  no. 1, 163--179.

%\bibitem{SV}
%J.A. Sanders and F. Verhulst,
%{\em Averaging Methods in Nonlinear Dynamical Systems},
%Springer-Verlag, New York 1985.

\bibitem{Sh94}
S. Schochet,
{\em Fast Singular Limits of Hyperbolic PDE's},
J. Diff. Equations {\bf 114} (1994), 476--512.

\bibitem{Te}
R. Temam,
{\em Navier-Stokes Equations: Theory and Numerical Analysis}.
AMS Chelsea, Providence, 2001.

\bibitem{Yong-1}
W-A, Yong,
{\em Entropy and global existence for hyperbolic balance laws}.
Arch. Ration. Mech. Anal.  {\bf } 172  (2004),  no. 2, 247--266.

\bibitem{Villani}
C. Villani,
{\em Hypocoercivity}, Memoirs Amer. Math. Soc. (to appear).

\end{thebibliography}
\end{document}